\documentclass[twoside, onecolumn, 10pt]{article}
\usepackage{latexsym}
\usepackage{amssymb}
\usepackage{mathrsfs}
\usepackage{amsmath}
\usepackage{marvosym}

\usepackage{pdfpages}

\newtheorem{theorem}{Theorem}[section]
\newtheorem{lemma}[theorem]{Lemma}
\newtheorem{corollary}[theorem]{Corollary}

\newtheorem{problem}[theorem]{Problem}

\begin{document}
\newcommand{\nb}{\nonumber}

\newcommand{\bb}{{\bf b}}
\newcommand{\bx}{{\bf x}}
\newcommand{\bv}{{\bf v}}

\def\bc{\begin{center}}
\def\ec{\end{center}}

\bc
{\LARGE
 Closed-form formulas for calculating the extremal ranks and inertias of a  quadratic matrix-valued function and their applications}
\ec

\begin{center}
{\large  Yongge Tian}
\end{center}

\begin{center}
{\em {\footnotesize CEMA, Central University of Finance and Economics, Beijing 100081, China
}}
\end{center}

\renewcommand{\thefootnote}{\fnsymbol{footnote}}
\footnotetext{ {\it E-mail Address:} yongge.tian@gmail.com}

\noindent {\small {\bf Abstract.} \
This paper presents a group of analytical formulas for calculating the global maximal
and minimal ranks and inertias of the quadratic matrix-valued function
$\phi(X) =  \left(\, AXB + C\,\right)\!M\!\left(\, AXB + C \right)^{*} + D$ and use them to derive
necessary and sufficient conditions for the two types of  multiple quadratic matrix-valued function
\begin{align*}
\left(\, \sum_{i = 1}^{k}A_iX_iB_i + C \,\right)\!M\!\left(\,\sum_{i = 1}^{k}A_iX_iB_i + C \, \right)^{\!*} +D,  \ \ \
\sum_{i = 1}^{k}\left(\,A_iX_iB_i + C_i\,\right)\!M_i\!\left(\,A_iX_iB_i  + C_i \,\right)^{*} +D
\end{align*}
to be semi-definite, respectively, where $A_i, \ B_i, \  C_i, \ C, \ D, \ M_i$
and $M$ are given matrices with $M_i$, $M$ and $D$ Hermitian, $i =1, \ldots, k$.  L\"owner partial ordering
optimizations of the two matrix-valued functions are studied and their solutions are characterized.

\medskip


\medskip

}

\section{Introduction}

\renewcommand{\theequation}{\thesection.\arabic{equation}}
\setcounter{section}{1}

This is the third part of the present author's work on quadratic matrix-valued functions and theii algebraic properties.
A matrix-valued function for complex matrices is a map between two matrix spaces ${\mathbb C}^{m \times n}$ and
${\mathbb C}^{p \times q}$, which can generally be written as
$$
Y = f(X) \ \  {\rm for}  \ \ Y \in {\mathbb C}^{m \times n} \ \ {\rm and} \ \  X \in {\mathbb C}^{p \times q},
$$
 or briefly, $f:  {\mathbb C}^{m \times n} \rightarrow {\mathbb C}^{p \times q}$. As usual, linear and quadratic matrix-valued
functions, as common representatives of various  matrix-valued functions, are extensively
studied from theoretical and applied points of view.

In this paper, we consider the following two types of
multiple quadratic matrix-valued function
\begin{align}
& \phi(\,X_1, \ldots, X_k\,) = \left(\, \sum_{i = 1}^{k}A_iX_iB_i + C \,\right)\!M\!\left(\,\sum_{i = 1}^{k}A_iX_iB_i + C \, \right)^{\!*} +D,
\label{ww11}
\\
& \psi(\,X_1, \ldots, X_k\,) = \sum_{i = 1}^{k}\left(\,A_iX_iB_i + C_i\,\right)\!M_i\!\left(\,A_iX_iB_i  + C_i \,\right)^{*} +D,
\label{ww12}
\end{align}
where  $A_i, \ B_i, \  C_i, \ C, \ D, \ M_i$ and $M$ are given matrices with $M_i$, $M$ and $D$ Hermitian, $X_i$ is a variable
matrix, $i =1, \ldots, k$. Eqs.\,(\ref{ww11}) or (\ref{ww12})  for $k =1$ is
\begin{align}
\phi(X)  = \left(\, AXB + C\,\right)\!M\!\left(\, AXB + C \right)^{*} + D,
\label{11}
\end{align}
where $A \in {\mathbb C}^{n\times p},$ $B \in {\mathbb C}^{m\times q},$  $C \in {\mathbb C}^{n\times q},$
$D \in {\mathbb C}_{{\rm H}}^{n}$ and $M \in {\mathbb C}_{{\rm H}}^{q}$  are given$,$ and $X\in {\mathbb C}^{p\times m}$
is a variable matrix$.$  We treat it as a combination $\phi = \tau\circ \rho$ of the following two simple linear and quadratic Hermitian matrix-valued functions:
\begin{align}
\rho
: X \rightarrow  AXB + C, \ \ \
\tau:  Y \rightarrow  YMY^{*} + D.
\label{11a}
\end{align}
For different choices of the given matrices, this quadratic function between matrix spaces includes many
ordinary quadratic forms and quadratic matrix-valued functions as its special cases, such as, ${\bf x}^{*}A{\bf x}$,
 $XAX^{*}$, $DXX^{*}D^{*}$, $(\, X - C\,)M(\,X - C)^{*}$, etc.

It is well known that quadratic functions in elementary mathematics and ordinary quadratic forms in linear algebra
 have a fairly complete theory with a long history and numerous applications.  Much of the beauty of
 these  quadratic objects were  highly appreciated by  mathematicians in all times, and
  and many of the fundamental ideas of quadratic functions and quadratic forms were developed
 in all branches of mathematics. While the mathematics of classic quadratic forms has been established
 for about one and a half century, various extensions of classic quadratic forms to some general settings were
 conducted from theoretical and applied point of view, in particular, quadratic matrix-valued functions and the corresponding quadratic matrix equations and quadratic matrix inequalities often appear briefly when needed to solve a variety of problems in mathematics and applications. These quadratic objects have many attractive features both from manipulative and computational point of view, and there is an intensive interest in studying behaviors of quadratic matrix-valued functions, quadratic matrix equations and quadratic matrix inequalities.  In fact, any essential development on the researches of quadratic objects will lead to many progresses in both mathematics and applications.  Compared with the theory of
ordinary  quadratic functions and forms, two distinctive features of quadratic matrix-valued functions
are the freedom of entries in variable matrices and the non-commutativity of matrix algebra.
 So that there is no a general theory for describing behaviors of a given quadratic matrix-valued function with multiple terms. In particular, to solve an optimization problem on a quadratic matrix-valued function is believed to be
 NP hard in general, and thus  there is a long way to go to establish a perfect theory on  quadratic
 matrix-valued functions. In recent years, Tian conducted a seminal study on quadratic matrix-valued functions
 in \cite{T-na,T-laa12}, which gave an initial quantitative understanding of the nature of  matrix rank and
 inertia optimization problems, in particular, a simple and precise linearization method was introduced for studying
 quadratic or nonlinear matrix-valued functions, and  many explicit formulas were established for calculating
the extremal ranks and inertias of some simple quadratic matrix-valued functions.
For applications of quadratic matrix-valued functions, quadratic matrix equations and quadratic matrix
inequalities in optimization theory, system and control theory, see the references given in \cite{T-na,T-laa12}.

Throughout this paper,
\begin{enumerate}
\item[]  ${\mathbb C}^{m\times n}$ stands for the set of all $m\times n$ complex
matrices;\

\item[] ${\mathbb C}_{{\rm H}}^{m}$
stands for the set of all $m\times m$ complex Hermitian matrices;

\item[] $A^{*}$,  $r(A)$ and ${\mathscr R}(A)$ stand
for the conjugate transpose, rank and range (column space) of a
matrix $A\in {\mathbb C}^{m\times n}$, respectively;

\item[]  $I_m$ denotes the identity matrix of order $m$;

\item[]  $[\, A, \, B\,]$ denotes a row
block matrix consisting of $A$ and $B$;

\item[]  the Moore--Penrose inverse of $A\in {\mathbb
C}^{m\times n}$, denoted by $A^{\dag}$, is defined to be the unique
solution $X$ satisfying the four matrix equations $AXA = A,$ $XAX =
X,$ $(AX)^{*} = AX$ and $(XA)^{*} =XA$;

\item[] the symbols $E_A$ and $F_A$ stand for $E_A = I_m - AA^{\dag}$ and $F_A = I_n-A^{\dag}A$;

\item[] an $X \in  {\mathbb C}^{n\times m}$ is called a $g$-inverse of
$A \in  {\mathbb C}^{m\times n}$, denoted by $A^{-}$, if it
satisfies  $AXA = A$;\\[-5mm]

\item[] an $X\in {\mathbb C}_{{\rm H}}^{m}$ is called a Hermitian $g$-inverse of
$A \in {\mathbb C}_{{\rm H}}^{m}$, denoted by $A^{\sim}$, if it satisfies $AXA = A$;
 called a reflexive Hermitian $g$-inverse of
$A \in {\mathbb C}_{{\rm H}}^{m}$, denoted by $A^{\sim}_r$, if it satisfies $AXA = A$ and
$XAX = X$;

\item[] $i_{+}(A)$ and $i_{-}(A)$, called the partial inertia of  $A \in {\mathbb C}_{{\rm H}}^{m}$,
 are defined to be  the numbers of the positive and negative eigenvalues of $A$ counted with
multiplicities, respectively;

$A \succ 0$ ($A \succcurlyeq  0, \, \prec 0, \, \preccurlyeq 0$) means that $A$ is Hermitian positive definite
(positive semi-definite, negative definite, negative semi-definite);

\item[]  two $A,  \, B \in {\mathbb C}_{{\rm H}}^{m}$ are said to satisfy the inequality $A \succ B$
($A \succcurlyeq B)$ in the L\"owner partial ordering if $A - B$ is  positive definite  (positive semi-definite).
\end{enumerate}

\section{Problem formulation}

\renewcommand{\theequation}{\thesection.\arabic{equation}}
\setcounter{section}{2}
\setcounter{equation}{0}

Matrix rank and inertia optimization problems are a class of discontinuous optimization problems,
in which decision variables are matrices running over certain matrix sets, while the rank and inertia
of the variable matrices are taken as integer-valued objective functions.
Because rank and inertia of matrices are always integers, no approximation methods are allowed to use
 when finding the maximal and minimal possible ranks and inertias of a matrix-valued function. So that
matrix rank and inertia optimization problems are not consistent with anyone of the ordinary continuous and
discrete problems in optimization theory. Less people paid attention to this kind of optimization problems,
 and no complete theory was established. But, the present author has been working on this topic with great
 effort in the past 30 years, and contribute a huge amount of results on matrix rank and inertia optimization
 problems.

 A major purpose of this paper is to develop a
unified optimization theory on ranks, inertias and partial orderings
 of quadratic matrix-valued functions by using pure algebraic operations of matrices,
 which enables us to handle many mathematical and applied problems on
 behaviors of quadratic matrix-valued functions, quadratic matrix equations and quadratic matrix inequalities.
The rank and inertia of a (Hermitian) matrix are two oldest basic concepts in linear algebra
for describing the dimension of the row/column vector space and the sign distribution
of the eigenvalues of the square matrix, which are well understood and are easy to compute by
the well-known elementary or congruent matrix operations. These two quantities
play an essential role in characterizing algebraic properties of (Hermitian) matrices.
Because concepts of  ranks and inertias are so generic in linear algebra, it is doubt that a primary
work in linear algebra is to establish (expansion) formulas for calculating ranks and inertias of matrices
as more as possible. However, this valuable work was really neglected in the development of linear algebra,
and a great chance for discovering thousands of rank and inertia formulas, some of which are given in Lemmas 3.2,
 3.3, 3.5 and 3.6 below,  were lost in the earlier period of linear algebra. This paper tries to make some
essential contributions on establishing formulas for ranks and inertias of some quadratic matrix-valued
functions.

Taking the rank and inertia of (\ref{11}) as inter-valued objective functions, we solve the following problems:

\begin{problem} \label{TS11}
{\rm
For the function in {\rm (\ref{11})},  establish explicit formulas for calculating
the following global extremal ranks and inertias
\begin{align}
 \max_{X\in {\mathbb C}^{p\times m}}r[\phi(X)], \ \ \  \min_{X\in {\mathbb C}^{p\times m}}r[\phi(X)], \ \ \
  \max_{X\in {\mathbb C}^{p\times m}}i_{\pm}[\phi(X)],
   \ \ \ \min_{X\in {\mathbb C}^{p\times m}}i_{\pm}[\phi(X)].
\label{12}
\end{align}
}
\end{problem}


\begin{problem} \label{T1S12}
{\rm
For the function in {\rm (\ref{11})}$,$
\begin{enumerate}
\item[{\rm (i)}]  establish necessary and sufficient conditions for the existence of an
$X \in {\mathbb C}^{p\times m}$ such that
\begin{align}
\phi(X)  =0;
\label{13}
\end{align}

\item[{\rm (ii)}] establish necessary and sufficient conditions for the following inequalities
\begin{align}
\phi(X) \succ 0, \ \  \phi(X) \succcurlyeq 0, \ \  \phi(X) \prec  0, \ \   \phi(X) \preccurlyeq 0
\label{14a}
\end{align}
to hold for an $X \in {\mathbb C}^{p\times m}$, respectively;

\item[{\rm (iii)}]  establish necessary and sufficient conditions for
\begin{align}
 \phi(X) \succ 0, \ \  \phi(X) \succcurlyeq 0, \ \  \phi(X) \prec  0, \ \   \phi(X) \preccurlyeq 0 \ \ \mbox{for all
  $X \in {\mathbb C}^{p\times m}$}
\label{14b}
\end{align}
to hold, respectively,  namely, to give identifying conditions for $\phi(X)$ to be a positive definite, positive semi-definite, negative definite, negative semi-definite function on complex matrices, respectively.
\end{enumerate}
}
\end{problem}

\begin{problem} \label{TS13}
{\rm For the function in {\rm (\ref{11})},  establish necessary and
sufficient conditions for the existence of  $\widehat{X}, \, \widetilde{X}
\in {\mathbb C}^{p\times m}$ such that
\begin{align}
\phi(X) \succcurlyeq \phi(\widehat{X}), \ \ \ \phi(X) \preccurlyeq \phi(\widetilde{X})
\label{13a}
\end{align}
hold for all $X\in {\mathbb C}^{p\times m}$, respectively, and derive analytical expressions
of the two matrices $\widehat{X}$ and $\widetilde{X}$. }
\end{problem}

\section{Preliminary results}

\renewcommand{\theequation}{\thesection.\arabic{equation}}
\setcounter{section}{3}
\setcounter{equation}{0}

The ranks and inertias are two generic indices in finite dimensional algebras. The results related  to
these indices are unreplaceable by any other quantitative tools in mathematics.  A simple but striking
 fact about the indices is stated in  following lemma.

\begin{lemma}  \label{T12}
Let ${\cal S}$ be a subset in ${\mathbb C}^{m\times n},$ and let ${\cal H}$ be a subset in
${\mathbb C}_{{\rm H}}^{m}.$ Then the following hold$.$
\begin{enumerate}
\item[{\rm (a)}]  Under $m =n,$ ${\cal S}$ has a nonsingular matrix if and only if
$\max_{X\in {\cal S}} r(X) = m.$

\item[{\rm (b)}] Under $m =n,$  all $X\in {\cal S}$ are nonsingular if and only if
 $\min_{X\in {\cal S}} r(X) = m.$

\item[{\rm (c)}] $0\in {\cal S}$ if and only if
$\min_{X\in {\cal S}} r(X) = 0.$

\item[{\rm (d)}] ${\cal S} = \{ 0\}$ if and only if $\max_{X\in {\cal S}} r(X) = 0.$

\item[{\rm (e)}] ${\cal H}$ has a matrix $X \succ 0$  $(X \prec 0)$ if and only if
$\max_{X\in {\cal H}} i_{+}(X) = m  \ \left(\, \max_{X\in {\cal H}} i_{-}(X) = m \,\right)\!.$

\item[{\rm (f)}] All $X\in {\cal H}$ satisfy $X \succ0$ $(X \prec 0),$ namely$,$  ${\cal H}$ is a subset of
the cone of positive definite matrices $($negative definite matrices$),$ if and only if
$\min_{X\in {\cal H}} i_{+}(X) = m \ \left(\,\min_{X\in {\cal H}} i_{-}(X) = m \, \right)\!.$

\item[{\rm (g)}] ${\cal H}$ has a matrix  $X \succcurlyeq 0$ $(X \preccurlyeq 0)$ if and only if
$\min_{X \in {\cal H}} i_{-}(X) = 0 \ \left(\,\min_{X \in {\cal H}} i_{+}(X) = 0 \,\right)\!.$

\item[{\rm (h)}] All $X\in {\cal H}$ satisfy $X \succcurlyeq 0$ $(X \preccurlyeq 0),$ namely$,$  ${\cal H}$ is a subset of
the cone of  positive semi-definite matrices $($\, negative semi-definite matrices\,$),$ if and only if
$\max_{X\in {\cal H}} i_{-}(X) = 0 \ \left(\, \max_{X\in {\cal H}} i_{+}(\,X) = 0 \,\right)\!.$
\end{enumerate}
\end{lemma}

The question of whether a given  matrix-valued function is semi-definite everywhere is ubiquitous in matrix theory and
applications. Lemma \ref{T12}(e)--(h) assert that if certain explicit formulas for calculating the global
maximal and minimal inertias of Hermitian matrix-valued functions are established, we can use them as a quantitative tool,
as demonstrated in Sections 2--7 below, to derive necessary and sufficient conditions for the
matrix-valued functions to be definite or semi-definite. In addition, we are able to use these inertia formulas to
establish various matrix inequalities in the L\"owner partial ordering, and to solve many
matrix optimization problems in the L\"owner partial ordering.

The following are obvious or well known (see \cite{T-na}), which will be used in the latter part of this paper.

\begin{lemma} \label{T13}
Let $A\in {\mathbb C}^m_{{\rm H}},$  $B \in {\mathbb C}^n_{{\rm
H}},$ $Q \in {\mathbb C}^{m \times n},$ and $P\in {\mathbb C}^{p\times m}$ with $r(P) = m.$ Then$,$
\begin{align}
& i_{\pm}(PAP^{*})  = i_{\pm}(A),
\label{15}
\\
& i_{\pm}(\lambda A)  = \left\{ \begin{array}{ll}  i_{\pm}(A)  &  if  \ \lambda  >0
 \\ i_{\mp}(A)  & if  \ \lambda < 0
 \end{array},
\right.
\label{16}
\\
i_{\pm}\!\left[ \begin{array}{cc}  A  & 0 \\  0  & B \end{array}
\right] & = i_{\pm}(A) + i_{\pm}(B),
\label{17}
\\
i_{+}\!\left[ \begin{array}{cc}  0  & Q \\  Q^{*} & 0 \end{array}
\right] & = i_{-}\!\left[ \begin{array}{cc}  0  & Q \\  Q^{*} & 0 \end{array}
\right] =r(Q).
\label{18}
\end{align}
\end{lemma}

\begin{lemma} [\cite{T-laa10}] \label{T14}
Let $A \in {\mathbb C}_{{\rm H}}^{m},$ $B \in \mathbb
C^{m\times n},$   $D \in {\mathbb C}_{{\rm H}}^{n},$   and let
$$
M_1 = \left[\!\!\begin{array}{cc}  A  & B  \\ B^{*}  & 0 \end{array}
\!\!\right]\!, \ \ M_2 = \left[\!\!\begin{array}{cc}  A  & B  \\
B^{*}  & D \end{array} \!\!\right]\!.
$$
Then$,$ the following expansion formulas hold
\begin{align}
& i_{\pm}(M_1) = r(B) + i_{\pm}(E_BAE_B),  \ \ \ \ \ \ \ \ \ \ \ \ \ \ \ \ \ \ \ \ r(M_1) = 2r(B) + r(E_BAE_B),
 \label{112k}
 \\
& i_{\pm}(M_2)  =i_{\pm}(A) + i_{\pm}\!\left[\!\!\begin{array}{cc} 0 &  E_AB
 \\
 B^{*}E_A & D - B^{*}A^{\dag}B \end{array}\!\!\right]\!, \ \
 r(M_2)  = r(A) + r\!\left[\!\!\begin{array}{cc} 0 &  E_AB
 \\
 B^{*}E_A & D - B^{*}A^{\dag}B \end{array}\!\!\right]\!.
\label{114k}
\end{align}
In particular$,$ the following hold.
\begin{enumerate}
\item[{\rm(a)}] If $A \succcurlyeq 0,$ then
\begin{align}
i_{+}(M_1)= r[\, A,  \, B \,], \ \ i_{-}(M_1) = r(B), \ \
 r(M_1) = r[\, A,  \, B \,] + r(B).
\label{116k}
\end{align}

\item[{\rm(b)}] If $A \preccurlyeq 0,$ then
\begin{align}
i_{+}(M_1) = r(B), \ \ i_{-}(M_1) = r[\, A,  \, B \,], \ \ r(M_1)
= r[\,A, \, B \,] + r(B).
\label{117k}
\end{align}

\item[{\rm(c)}] If ${\mathscr R}(B) \subseteq {\mathscr R}(A),$ then
\begin{align}
i_{\pm}(M_2) & = i_{\pm}(A) + i_{\pm}(\,  D - B^{*}A^{\dag}B \,), \ \
r(M_2) = r(A) + r(\,  D - B^{*}A^{\dag}B \,).
\label{118k}
\end{align}

\item[{\rm(d)}] $r(M_2) = r(A) \Leftrightarrow {\mathscr R}(B) \subseteq {\mathscr R}(A)$ and
$D = B^{*}A^{\dag}B.$

\item[{\rm(e)}] $M_2  \succcurlyeq 0 \Leftrightarrow A  \succcurlyeq 0, \   {\mathscr R}(B) \subseteq {\mathscr R}(A)$
 and $D  - B^{*}A^{\dag}B  \succcurlyeq 0.$
\end{enumerate}
\end{lemma}

\begin{lemma} [\cite{Pen}] \label{T15}
Let $A \in \mathbb C^{m \times p}, \, B \in \mathbb C^{q \times n}$
 and $C \in \mathbb C^{m \times n}$ be given$.$  Then the matrix equation $AXB = C$ is consistent if and only if
 ${\mathscr R}(C) \subseteq {\mathscr R}(A)$ and ${\mathscr R}(C^{*}) \subseteq {\mathscr R}(B^{*}),$
or equivalently$,$  $AA^{\dag}CB^{\dag}B = C.$  In this case$,$ the general solution can be written as
\begin{align}
X= A^{\dag}CB^{\dag} + F_AV_1 + V_2E_B,
\label{115a}
\end{align}
where $V_1$ and $V_2$ are arbitrary matrices$.$ In particular$,$ $AXB = C$ has a unique solution if and only if
\begin{align}
r(A) = p, \  \ r(B) = q, \ \ {\mathscr R}(C) \subseteq {\mathscr R}(A), \ \ {\mathscr R}(C^{*}) \subseteq {\mathscr R}(B^{*}).
\label{uu123}
\end{align}
\end{lemma}

\begin{lemma} [\cite{DG,Ti-se,TC}] \label{T15a}
Let $A \in \mathbb C^{m \times n}, \, B \in \mathbb C^{m \times p}$
 and $C \in \mathbb C^{q \times n}$ be given$,$ and  $X \in \mathbb C^{p\times q}$ be
 a variable matrix$.$  Then$,$  the global maximal and minimal ranks of  $A + BXC$ are given by
\begin{align}
\max_{X \in \mathbb C^{p \times q}}\!\!\!r(\, A  + BXC \,) & =  \min  \left\{ r[ \, A, \,  B \, ], \ \
 r \!\left[\!\! \begin{array}{c}  A  \\ C  \end{array}\!\!\right] \right\}\!,
\label{vv115}
\\
\min_{X \in \mathbb C^{p \times q}} \!\!\!r(\, A + BXC \,)  & =  r[\, A, \, B \,] +
 r \!\left[\!\! \begin{array}{c} A \\ C \end{array} \!\!\right] -
 r \!\left[\!\! \begin{array}{cc} A & B \\ C & 0 \end{array} \!\!\right]\!.
\label{vv116}
\end{align}
\end{lemma}

\begin{lemma}  [\cite{LT-jota}] \label{T16}
Let  $A \in {\mathbb C}_{{\rm H}}^{m}$, $B \in \mathbb
C^{m\times n}$ and $C \in \mathbb C^{p \times m}$ be given$,$
 $X \in \mathbb C^{n\times p}$ be a variable matrix$.$
  Then$,$
\begin{align}
&\max_{X\in \mathbb C^{n \times p}}\!\!\!r[\,  A  + BXC  + (BXC)^{*}\,]
 = \min \left\{ r[\,A,\, B,\, C^{*}\,],
 \ \ r\!\left[\begin{array}{cc} A & B
\\B^{*} & 0
\end{array}\right]\!, \ \ r\!\left[\!\begin{array}{cc} A & C^{*}
\\ C & 0
\end{array}\!\right] \right\}\!,
\label{120}
\\
&\min_{X\in \mathbb C^{n \times p}}\!\!\!r[\, A  + BXC  + (BXC)^{*} \,]  =
2r[\,A,\, B,\, C^{*}\,] + \max\{\, s_{+} + s_{-}, \
t_{+} + t_{-}, \ s_{+} + t_{-}, \ s_{-} + t_{+}  \},
\label{121}
\\
& \max_{X\in \mathbb C^{n \times p}}\!\!\!i_{\pm}[\,  A  + BXC  + (BXC)^{*}\,]
= \min\!\left\{i_{\pm}\!\left[\!\begin{array}{ccc} A
 & B  \\  B^{*}  & 0
   \end{array}\!\right]\!, \ \ i_{\pm}\!\left[\!\begin{array}{ccc} A
 & C^{*}  \\  C  & 0
   \end{array}\!\right] \right\}\!,
\label{122}
\\
& \min_{X\in \mathbb C^{n \times p}}\!\!\!i_{\pm}[\,  A  + BXC  + (BXC)^{*}\,]
= r[\,A,\, B,\, C^{*}\,] + \max\{\, s_{\pm}, \ \ t_{\pm} \, \},
\label{123}
\end{align}
where
\begin{align*}
s_{\pm}  = i_{\pm}\!\left[\!\!\begin{array}{cc} A & B \\ B^{*} & 0\end{array}\!\!\right]
 - r\!\left[\begin{array}{ccc} A & B  & C^{*} \\ B^{*} & 0 & 0
\end{array}\!\!\right]\!, \ \
 t_{\pm}  =i_{\pm}\!\left[\!\!\begin{array}{cc} A & C^{*} \\ C & 0\end{array}\!\!\right]
 - r\!\left[\begin{array}{ccc} A & B  & C^{*} \\ C & 0 & 0 \end{array}\!\!\right]\!.
\end{align*}

\end{lemma}

\section{Main results}
\renewcommand{\theequation}{\thesection.\arabic{equation}}
\setcounter{section}{4}
\setcounter{equation}{0}


We first solve Problem 1.1 through a linearization method and Theorem 1.8.

\begin{theorem}\label{T41}
Let $ \phi(X)$  be as given  in {\rm (\ref{11})}$,$  and define
\begin{align}
& N_1 =\left[\!\! \begin{array}{ccc} D  + CMC^{*} & A \\ A^{*} & 0
 \end{array} \!\!\right]\!,  \ \ N_2 = \left[\!\! \begin{array}{cccc} D  + CMC^{*} & CMB^{*}  & A
 \\ A^{*} & 0 & 0 \end{array} \!\!\right]\!,
\label{21}
 \\
&   N_3 = \left[\!\! \begin{array}{ccc} D  + CMC^{*} & CMB^{*} \\
 BMC^{*} & BMB^{*} \end{array} \!\!\right]\!, \ \ N_4 = \left[\!\! \begin{array}{cccc}
  D  + CMC^{*} & CMB^{*}  & A \\ BMC^{*} & BMB^{*} & 0  \end{array} \!\!\right]\!.
\label{22}
\end{align}
Then$,$ the global maximal  and minimal ranks and inertias of $\phi(X)$ are given by
\begin{align}
\max_{X\in {\mathbb C}^{p \times m}} \!\!\!r[\,\phi(X)\,]  & =  \min  \left\{ \, r[\, D  + CMC^{*}, \, CMB^{*}, \, A \,], \ \
r(N_1),  \ \ r(N_3) \, \right\},
\label{23}
\\
\min_{X\in {\mathbb C}^{p \times m}} \!\!\!r[\,\phi(X)\,] & = 2r[\, D  + CMC^{*}, \, CMB^{*}, \, A \,]  +
\max\{\,  s_1, \ \  s_2,  \ \ s_3, \ \  s_4 \,\},
\label{24}
\\
\max_{X\in {\mathbb C}^{p \times m}} \!\!\!i_{\pm}[\,\phi(X)\,] &  =
 \min \left\{ \, i_{\pm}(N_1),  \ \ i_{\pm}(N_3) \, \right\},
\label{25}
\\
\min_{X\in {\mathbb C}^{p \times m}} \!\!\!i_{\pm}[\,\phi(X)\,] & = r[\, D  + CMC^{*}, \, CMB^{*}, \, A \,]  +
\max\left\{\,  i_{\pm}(N_1) - r(N_2), \ \  i_{\pm}(N_3) - r(N_4) \, \right\},
 \label{26}
\end{align}
where
\begin{align*}
& s_1 = r(N_1) - 2r(N_2), \ \  s_2 = r(N_3) - 2r(N_4),
\\
& s_3 = i_{+}(N_1) + i_{-}(N_3) - r(N_2) - r(N_4), \ \  s_4 = i_{-}(N_1) + i_{+}(N_3) - r(N_2) - r(N_4).
\end{align*}
\end{theorem}

\noindent  {\bf Proof.} \
It is easy to verify from (\ref{114k}) that
\begin{align}
& i_{\pm}[\,\left(\, AXB + C\,\right)\!M\!\left(\, AXB + C \right)^{*} + D\,]
= i_{\pm}\!\left[\!\! \begin{array}{cc} -M & M(\, AXB + C \,)^{*} \\
(\, AXB + C \,)M  & D  \end{array} \!\!\right]  -  i_{\pm}(-M),
\label{27}
\\
& r[\, \left(\, AXB + C\,\right)\!M\!\left(\, AXB + C \right)^{*} + D\,]
= r\!\left[\!\! \begin{array}{cc} -M & M(\, AXB + C \,)^{*} \\
(\, AXB + C \,)M  & D  \end{array} \!\!\right]   -  r(M).
\label{28}
\end{align}
that is, the rank and inertia of $\phi(X)$ in  (\ref{11}) can be calculated by those of the following
linear matrix-valued function
\begin{align}
\psi(X) & = \left[\!\! \begin{array}{cc} -M & M\!\left(\, AXB + C \right)^{*} \\
(\, AXB + C \,)M  & D  \end{array} \!\!\right]  \nb
\\
&  = \left[\!\! \begin{array}{cc} -M & MC^{*}
\\
CM  & D  \end{array} \!\!\right]  + \left[\!\! \begin{array}{c} 0 \\ A \end{array} \!\!\right]X[\, BM, \, 0\,]
 + \left[\!\! \begin{array}{c} MB^{*} \\ 0 \end{array} \!\!\right]X^{*}[\, 0, \, A^{*} \,].
\label{29}
\end{align}
Note from (\ref{27}) and (\ref{28}) that
\begin{align}
& \max_{X\in {\mathbb C}^{p \times m}} \!\!\!r[\,\phi(X)\,] =
\max_{X\in {\mathbb C}^{p \times m}}\!\!r[\psi(X)]  -  r(A),
\label{210}
\\
&  \min_{X\in {\mathbb C}^{p \times m}} \!\!\!r[\,\phi(X)\,]
= \min_{X\in {\mathbb C}^{p \times m}}\!\!r[\psi(X)]  -  r(A),
\label{211}
\\
& \max_{X\in {\mathbb C}^{p \times m}} \!\!\!i_{\pm}[\,\phi(X)\,]
= \max_{X\in {\mathbb C}^{p \times m}}\!\!i_{\pm}[\psi(X)]
-  i_{\mp}(A),
\label{212}
\\
& \min_{X\in {\mathbb C}^{p \times m}} \!\!\!i_{\pm}[\,\phi(X)\,] =
\min_{X\in {\mathbb C}^{p \times m}}\!\!\!i_{\pm}[\psi(X)]  -  i_{\mp}(A).
\label{213}
\end{align}
Applying Lemma \ref{T16} to (\ref{29}), we first obtain
\begin{align}
\max_{X\in {\mathbb C}^{p \times m}}\!\!\!r[\psi(X)] & =\min\{\, r(H),\ \ r(G_1), \ \ r(G_2) \,\},
\label{214}
\\
\min_{X\in {\mathbb C}^{p \times m}}\!\!\!r[\psi(X)] & =2r(H)+\max\{\, s_{+} + s_{-}, \ \
t_{+} + t_{-},  \ \ s_{+} + t_{-}, \ \ s_{-} + t_{+} \},
\label{215}
\\
\max_{X\in {\mathbb C}^{p \times m}}\!\!i_{\pm} [\psi(X)] & = \min\{\, i_{\pm}(G_1), \ \ i_{\pm}(G_2) \, \},
\label{216}
\\
 \min_{X\in {\mathbb C}^{p \times m}} \!\!i_{\pm}[\psi(X)] &  = r(H) + \max\{\, s_{\pm},\ \ t_{\pm} \, \},
\label{217}
\end{align}
where
$$
H= \left[\!\! \begin{array}{cccc} -M & MC^{*} & 0 & MB^{*} \\ CM  & D & A & 0 \end{array} \!\!\right]\!,
 \  G_1 = \!\left[\!\! \begin{array}{ccc}  -M  & MC^{*}  & 0  \\ CM  & D & A \\  0 & A^{*} & 0
 \end{array} \!\!\right], \
 G_2 = \left[\!\! \begin{array}{cccc}   -M  & MC^{*}  & MB^{*} \\ CM  & D & 0 \\  BM & 0 & 0 \end{array}
\!\!\right]\!,
$$\\[-6mm]
$$
H_1 = \left[\!\! \begin{array}{cccc} -M  & MC^{*}  & 0 & MB^{*}  \\ CM  & D & A & 0\\  0 & A^{*} & 0 &0 \end{array} \!\!\right]\!, \
H_2 = \left[\!\! \begin{array}{cccc}   -M  & MC^{*}  & MB^{*}  & 0 \\ CM  & D & 0 & A \\  BM & 0 & 0 &0 \end{array} \!\!\right]\!,\
$$
and
$$
s_\pm = i_\pm(G_1) - r(H_1), \ \  t_\pm =i_\pm(G_2) - r(H_2).
$$
 It is easy to derive from Lemmas \ref{T13} and \ref{T14},
 elementary  matrix operations and congruence matrix operations  that
\begin{align}
& r(H)  =  r(M) + r[\, D  + CMC^{*}, \, CMB^{*}, \, A \,],
\label{218}
\\
&  r(H_1) = r\!\left[\!\! \begin{array}{cccc} D  + CMC^{*} & CMB^{*}  & A
 \\ A^{*} & 0 & 0 \end{array} \!\!\right]\!, \  r(H_2) = r(M) +  r\!\left[\!\! \begin{array}{cccc}
  D  + CMC^{*} & CMB^{*}  & A \\ BMC^{*} & BMB^{*} & 0  \end{array} \!\!\right]\!,
\label{219}
 \\
&  i_\pm(G_1)= i_{\mp}(M) + i_{\pm}\!\left[\!\! \begin{array}{ccc} C + CMC^{*} & A \\ A^{*} & 0
 \end{array} \!\!\right]\!,  \ i_\pm(G_2)= i_{\mp}(M) +
 i_{\pm}\!\left[\!\! \begin{array}{ccc} D  + CMC^{*} & CMB^{*} \\
 BMC^{*} & BMB^{*} \end{array} \!\!\right]\!.
\label{220}
\end{align}
Hence,
\begin{align}
& r(G_1)= r(M) + r(N_1), \ \ r(G_2)= r(M) + r(N_3),
\label{221}
\\
& s_{\pm}  = i_\pm (G_1) - r(H_1) = i_{\pm}(N_1) - r(N_2) - i_{\pm}(M),
\label{222}
\\
& t_{\pm}  = i_\pm (G_2) - r(H_2)= i_{\pm}(N_3) - r(N_4) -  i_{\pm}(M).
\label{223}
\end{align}
Substituting (\ref{218})--(\ref{223}) into (\ref{215})--(\ref{217}), and then
(\ref{214})--(\ref{217}) into  (\ref{210})--(\ref{213}), we obtain (\ref{24})--(\ref{27}). \qquad $\Box$

\medskip

Without loss of generality, we assume in what follows that both $A \neq 0$ and $BMB^{*} \neq 0$ in (\ref{11}).
Applying Lemma 1.4 to (\ref{23})--(\ref{26}), we obtain the following results.

\begin{corollary}\label{T42}
Let $\phi(X)$ be as given  in {\rm (\ref{11})}$,$ and $N_1$ and $N_3$ be the matrices of {\rm (\ref{21})} and
 {\rm (\ref{22})}$.$  Then$,$ the following hold$.$
\begin{enumerate}
\item[{\rm (a)}]  There exists an $X \in {\mathbb C}^{p \times m}$ such that $\phi(X)$ is
 nonsingular if and only if $r[\, D + CMC^{*}, \, CMB^{*}, \, A \,] =n,$  $r(N_1) \geqslant n$ and $r(N_3) \geqslant n.$

\item[{\rm (b)}] $\phi(X)$ is nonsingular for all $X \in {\mathbb C}^{p \times m}$ if and only if
$r(\, D + CMC^{*} \,) = n,$  and one of the following four conditions holds
\begin{enumerate}
\item[{\rm (i)}]  $BMC^{*}(\,D + CMC^{*}\,)^{-1}A =0$ and $A^{*}(\,D + CMC^{*}\,)^{-1}A =0.$

\item[{\rm (ii)}] $BMC^{*}(\,D + CMC^{*}\,)^{-1}A =0$ and $BMC^{*}(\,D + CMC^{*}\,)^{-1}CMB^{*} =BMB^{*}.$

\item[{\rm (iii)}] $A^{*}(\,D + CMC^{*}\,)^{-1}A \succcurlyeq 0, \  BMB^{*} - BMC^{*}(\,D + CMC^{*}\,)^{-1}CMB^{*} \succcurlyeq 0, \
{\mathscr R}[\, A^{*}(\,D + CMC^{*}\,)^{-1}CMB^{*}] \subseteq  {\mathscr R}[\, A^{*}(\,D + CMC^{*}\,)^{-1}A\,]$ and
${\mathscr R}[\,BMC^{*}(\,D + CMC^{*}\,)^{-1}A\,] \subseteq  {\mathscr R}[\, BMB^{*} - BMC^{*}(\,D + CMC^{*}\,)^{-1}CMB^{*}\,].$

 \item[{\rm (iv)}]  $A^{*}(\,D + CMC^{*}\,)^{-1}A \preccurlyeq 0, \  BMB^{*} - BMC^{*}(\,D + CMC^{*}\,)^{-1}CMB^{*} \preccurlyeq 0, \
{\mathscr R}[\,A^{*}(\,D + CMC^{*}\,)^{-1}CMB^{*}] \subseteq  {\mathscr R}[\,A^{*}(\,D + CMC^{*}\,)^{-1}A\,],$ and
${\mathscr R}[\,BMC^{*}(\,D + CMC^{*}\,)^{-1}A\,] \subseteq  {\mathscr R}[\, BMB^{*} - BMC^{*}(\,D + CMC^{*}\,)^{-1}CMB^{*}\,].$
\end{enumerate}

\item[{\rm (c)}] There exists an $X \in {\mathbb C}^{p \times m}$ such that $\phi(X) =0,$  namely$,$ the
matrix equation in {\rm (\ref{13})} is consistent$,$ if and only if
\begin{align*}
& {\mathscr R}(\,D + CMC^{*}\,) \subseteq {\mathscr R}[\, A, \, CMB^{*} \,], \ \ \  r(M_1) = 2r(A), \ \ \
2r[\,A, \,  CMB^{*}\,] + r(N_3) -2r(N_4) \preccurlyeq  0,  \ \ \
\\
&  r[\,A, \,  CMB^{*}\,] + i_{+}(N_3) - r(N_4)  \leqslant 0,
\ \ \ r[\,A, \,  CMB^{*}\,] +   i_{-}(N_3) - r(N_4) \leqslant 0.
\end{align*}

\item[{\rm (d)}] There exists an $X  \in {\mathbb C}^{p \times m}$ such that $\phi(X)\succ0,$ namely$,$
the matrix inequality is feasible$,$ if and only if
$$
i_{+}(N_1) = n \ and \ i_{+}(N_3) \geqslant n, \ \ or  \ \  i_{+}(N_1) \geqslant  n \ and \
i_{+}(N_3) = n.
$$

\item[{\rm (e)}] There exists an $X  \in {\mathbb C}^{p \times m}$ such that $\phi(X)\prec0,$
the matrix inequality is feasible$,$ if and only if
$$
i_{-}(N_1) = n  \ and \ i_{-}(N_3) \geqslant n, \ \   or \ \ i_{-}(N_1) \geqslant n \ and \
i_{-}(N_3) = n.
$$

\item[{\rm (f)}] $\phi(X)\succ0$ for all $X \in {\mathbb C}^{p \times m}$, namely$,$ $\phi(X)$ is a completely
positive matrix-valued function$,$
if and only if
\begin{align*}
D + CMC^{*} \succ 0,  \  N_3 \succcurlyeq 0, \  {\mathscr R}\left[\!\! \begin{array}{c} A \\ 0
\end{array} \!\!\right]  \subseteq  {\mathscr R}(N_3).
 \end{align*}

 \item[{\rm (g)}] $\phi(X)\prec0$ for all $X \in {\mathbb C}^{p \times m}$ namely$,$ $\phi(X)$ is
  a completely negative  matrix-valued function$,$
 if and only if
$$
D + CMC^{*} \prec 0,  \  N_3 \preccurlyeq 0, \ {\mathscr R}\left[\!\! \begin{array}{c} A \\ 0 \end{array} \!\!\right]  \subseteq
 {\mathscr R}(N_3).
$$

\item[{\rm (h)}] There exists an $X  \in {\mathbb C}^{p \times m}$ such that $\phi(X)\succcurlyeq 0,$ namely$,$
the matrix inequality is feasible$,$ if and only if
$$
r[\, D + CMC^{*}, \, CMB^{*}, \, A \,] + i_{-}(N_1)\leqslant r(N_2) \  and \  r[\, D + CMC^{*}, \, CMB^{*}, \, A \,] +
i_{-}(N_3)\leqslant r(N_4).
$$

\item[{\rm (i)}] There exists an $X  \in {\mathbb C}^{p \times m}$ such that $\phi(X)\preccurlyeq 0,$ namely$,$
the matrix inequality is feasible$,$ if and only if
$$
r[\, D + CMC^{*}, \, CMB^{*}, \, A \,] + i_{+}(N_1)\leqslant r(N_2) \  and \  r[\, D + CMC^{*}, \, CMB^{*}, \, A \,] +
i_{+}(N_3)\leqslant r(N_4).
$$

\item[{\rm (j)}] $\phi(X)\succcurlyeq 0$ for all $X \in {\mathbb C}^{p \times m},$
namely$,$ $\phi(X)$ is  a positive a positive matrix-valued  function$,$  if and only if $N_3  \succcurlyeq 0.$

  \item[{\rm (k)}] $\phi(X)\preccurlyeq 0$ for all $X \in {\mathbb C}^{p \times m},$
namely$,$ $\phi(X)$ is  a negative matrix-valued function$,$  if and only if $N_3 \preccurlyeq 0.$
 \end{enumerate}
\end{corollary}

\noindent  {\bf Proof.} \ We only show (b). Under the condition $r(\,D + CMC^{*}\,) = n,$   (\ref{24}) reduces to
\begin{align}
\min_{X\in {\mathbb C}^{p \times m}} \!\!\!r[\,\phi(X)\,] = 2n +
\max\{\,  s_1, \ \  s_2,  \ \ s_3, \ \  s_4 \,\},
\label{224}
\end{align}
where
\begin{align*}
& s_1 =  r[A^{*}(\,D + CMC^{*}\,)^{-1}A] -  2r[\, A^{*}(\,D + CMC^{*}\,)^{-1}CMB^{*}, \, A^{*}(\,D + CMC^{*}\,)^{-1}A \,] -n,
\\
& s_2 =  r[\, BMB^{*} - BMC^{*}(\,D + CMC^{*}\,)^{-1}CMB^{*}\,]
\\
& \ \ \ \ \ \   - 2r[\, BMB^{*}- BMC^{*}(\,D + CMC^{*}\,)^{-1}CMB^{*}, \, BMC^{*}(\,D + CMC^{*}\,)^{-1}A \,] -n,
\\
& s_3 = i_{-}[A^{*}(\,D + CMC^{*}\,)^{-1}A] + i_{-}[\, BMB^{*} - BMC^{*}(\,D + CMC^{*}\,)^{-1}CMB^{*}\,]
\\
&  \ \ \ \ - r[\, A^{*}(\,D + CMC^{*}\,)^{-1}CMB^{*}, \, A^{*}(\,D + CMC^{*}\,)^{-1}A \,]
\\
&  \ \ \ \ - r[\, BMB^{*} - BMC^{*}(\,D + CMC^{*}\,)^{-1}CMB^{*}, \, BMC^{*}(\,D + CMC^{*}\,)^{-1}A \,] -n,
\\
& s_4 =  i_{+}[A^{*}(\,D + CMC^{*}\,)^{-1}A] + i_{+}[\, BMB^{*} - BMC^{*}(\,D + CMC^{*}\,)^{-1}CMB^{*}\,]
 \\
 & \ \ \ \  - r[\, A^{*}(\,D + CMC^{*}\,)^{-1}CMB^{*}, \, A^{*}(\,D + CMC^{*}\,)^{-1}A \,]
\\
& \ \ \ \ - r[\, BMB^{*} - BMC^{*}(\,D + CMC^{*}\,)^{-1}CMB^{*}, \, BMC^{*}(\,D + CMC^{*}\,)^{-1}A \,] -n.
\end{align*}
Setting (\ref{224}) equal to $n$, we see that
$\phi(X)$ is nonsingular for all $X \in {\mathbb C}^{p \times m}$ if and only if
$r(D + CMC^{*}) = n,$  and one of the following four rank equalities holds
\begin{enumerate}
\item[{\rm (i)}] $r[A^{*}(\,D + CMC^{*}\,)^{-1}A] = 2r[\, A^{*}(\,D + CMC^{*}\,)^{-1}CMB^{*}, \, A^{*}(\,D + CMC^{*}\,)^{-1}A \,];$

\item[{\rm (ii)}] $ r[\, BMB^{*} - BMC^{*}(\,D + CMC^{*}\,)^{-1}CMB^{*}\,]$

$ = 2r[\, BMB^{*}- BMC^{*}(\,D + CMC^{*}\,)^{-1}CMB^{*}, \, BMC^{*}(\,D + CMC^{*}\,)^{-1}A \,];$

\item[{\rm (iii)}]  $i_{-}[A^{*}(\,D + CMC^{*}\,)^{-1}A] + i_{-}[\, BMB^{*} - BMC^{*}(\,D + CMC^{*}\,)^{-1}CMB^{*}\,]$ \\
$=r[\, A^{*}(\,D + CMC^{*}\,)^{-1}CMB^{*}, \, A^{*}(\,D + CMC^{*}\,)^{-1}A \,]$
\\
$ \ \ \ \ + \ r[\, BMB^{*} - BMC^{*}(\,D + CMC^{*}\,)^{-1}CMB^{*}, \, BMC^{*}(\,D + CMC^{*}\,)^{-1}A \,];$

\item[{\rm (iv)}] $i_{+}[A^{*}(\,D + CMC^{*}\,)^{-1}A] + i_{+}[\, BMB^{*} - BMC^{*}(\,D + CMC^{*}\,)^{-1}CMB^{*}\,]$

$=r[\, A^{*}(\,D + CMC^{*}\,)^{-1}CMB^{*}, \, A^{*}(\,D + CMC^{*}\,)^{-1}A \,]$

$ \ + \ r[\, BMB^{*} - BMC^{*}(\,D + CMC^{*}\,)^{-1}CMB^{*}, \, BMC^{*}(\,D + CMC^{*}\,)^{-1}A \,] -n;$
\end{enumerate}
which are further equivalent to the result in (b) by comparing both sides of the four rank equalities. \qquad $\Box$

\medskip

A special case of (\ref{11}) is
\begin{align}
\phi(X)  = \left(\, AXB + C\,\right)\left(\, AXB + C \right)^{*} - I_n,
\label{vv25}
\end{align}
where $\phi(X)  =0$ means that the rows of $AXB + C$ are orthogonal. Further, if
 $AXB + C$ is square, $\phi(X)  =0$ means that $AXB + C$ is unitary.
Applying Theorem 4.1 and Corollary 4.2 to (\ref{vv25}) will yield a group of consequences.

\begin{theorem}\label{T43nn}
Let $ \phi(X)$  be as given  in {\rm (\ref{11})}$,$  and define
\begin{align}
& N_1 =\left[\!\! \begin{array}{ccc} CC^{*} - I_n & A \\ A^{*} & 0
 \end{array} \!\!\right]\!,  \ \ N_2 = \left[\!\! \begin{array}{cccc}  CC^{*} - I_n & CB^{*}  & A
 \\ A^{*} & 0 & 0 \end{array} \!\!\right]\!,
\label{21aa}
 \\
&   N_3 = \left[\!\! \begin{array}{ccc}  CC^{*} - I_n & CB^{*} \\
 BC^{*} & BB^{*} \end{array} \!\!\right]\!, \ \ N_4 = \left[\!\! \begin{array}{cccc}
   CC^{*} - I_n & CB^{*}  & A \\ BC^{*} & BB^{*} & 0  \end{array} \!\!\right]\!.
\label{22aa}
\end{align}
Then$,$ the global maximal  and minimal ranks and inertias of $\phi(X)$ are given by
\begin{align}
\max_{X\in {\mathbb C}^{p \times m}} \!\!\!r[\,\phi(X)\,]  & =  \min  \left\{ \, r[\, CC^{*} -I_n, \, CB^{*}, \, A \,], \ \
r(N_1),  \ \ r(N_3) \, \right\},
\label{23aa}
\\
\min_{X\in {\mathbb C}^{p \times m}} \!\!\!r[\,\phi(X)\,] & = 2r[\,  CC^{*} - I_n, \, CB^{*}, \, A \,]  +
\max\{\,  s_1, \ \  s_2,  \ \ s_3, \ \  s_4 \,\},
\label{24aa}
\\
\max_{X\in {\mathbb C}^{p \times m}} \!\!\!i_{\pm}[\,\phi(X)\,] &  =
 \min \left\{ \, i_{\pm}(N_1),  \ \ i_{\pm}(N_3) \, \right\},
\label{25aa}
\\
\min_{X\in {\mathbb C}^{p \times m}} \!\!\!i_{\pm}[\,\phi(X)\,] & = r[\, CC^{*} -I_n, \, CB^{*}, \, A \,]  +
\max\left\{\,  i_{\pm}(N_1) - r(N_2), \ \  i_{\pm}(N_3) - r(N_4) \, \right\},
 \label{26aa}
\end{align}
where
\begin{align*}
& s_1 = r(N_1) - 2r(N_2), \ \  s_2 = r(N_3) - 2r(N_4),
\\
& s_3 = i_{+}(N_1) + i_{-}(N_3) - r(N_2) - r(N_4), \ \  s_4 = i_{-}(N_1) + i_{+}(N_3) - r(N_2) - r(N_4).
\end{align*}
\end{theorem}

Whether a given function is positive or nonnegative everywhere is a fundamental research subject in both
elementary and advanced mathematics. It was realized in matrix theory that the complexity status of the definite and semi-definite feasibility problems of a general matrix-valued function is NP-hard. Corollary \ref{T42}(d)--(k),
however, show that we are really able to characterize the definiteness and semi-definiteness of
(\ref{11}) by using some ordinary and elementary
methods. These results set up a criterion for characterizing definiteness and semi-definiteness of nonlinear
matrix-valued functions, and will prompt more investigations on this challenging topic. In particular,  definiteness and semi-definiteness of some  nonlinear matrix-valued functions generated from  (\ref{11}) can be identified. We shall present them in another paper.

Recall that a Hermitian matrix $A$ can uniquely be decomposed as the difference of two disjoint Hermitian positive
semi-definite definite  matrices
\begin{align}
A = A_1 - A_2,  \ \ A_1A_2 = A_2A_1 =0, \ \ A_1\succcurlyeq 0,  \ \ A_2\succcurlyeq 0 .
\label{zx225}
\end{align}

Applying this assertion to (\ref{11}),  we obtain the following result.

\begin{corollary}\label{T43}
Let $\phi(X)$ be as given  in {\rm (\ref{11})}$.$ Then$,$ $\phi(X)$ can always be decomposed as
\begin{align}
\phi(X) = \phi_1(X) - \phi_2(X),
\label{zx226}
\end{align}
where
$$
\phi_1(X)  = \left(\, AXB + C\,\right)\!M_1\!\left(\, AXB + C \right)^{*} + D_1, \ \
\phi_2(X)  = \left(\, AXB + C\,\right)\!M_2\!\left(\, AXB + C \right)^{*} + D_2
$$
satisfy
\begin{align}
\phi_1(X) \succcurlyeq 0, \ \  \phi_2(X) \succcurlyeq 0
\label{zx227}
\end{align}
 for all $X \in {\mathbb C}^{p \times m}.$
\end{corollary}

\noindent  {\bf Proof.} \  Note from (\ref{zx225}) that the two Hermitian matrices $D$ and $M$ in    (\ref{11})
can uniquely be  decomposed as
\begin{align*}
& D = D_1 - D_2,  \ \ D_1D_2 = D_2D_1 =0, \ \ D_1\succcurlyeq 0,  \ \ D_2\succcurlyeq 0,
\\
& M = M_1 - M_2,  \ \ M_1M_2 = M_2M_1 =0, \ \ M_1\succcurlyeq 0,  \ \ M_2\succcurlyeq 0,
\end{align*}
So that $\phi_1(X)$ and $\phi_2(X)$  in(\ref{zx226}) are positive matrix-valued functions.  \qquad $\Box$

\medskip

\begin{corollary}\label{T43a}
Let $\phi(X)$ be as given in {\rm (\ref{11})}$,$  and suppose that $AXB + C =0$ has a solution$,$
i.e.$,$ ${\mathscr R}(C) \subseteq {\mathscr R}(A)$ and ${\mathscr R}(C^{*}) \subseteq {\mathscr R}(B^{*}),$
and let $N =\left[\!\! \begin{array}{ccc} D & A \\ A^{*} & 0 \end{array} \!\!\right]\!.$ Then$,$
\begin{align}
\max_{X\in {\mathbb C}^{p \times m}}\!\!\! r[\, \phi(X)\,] & =  \min  \left\{ \, r[\, A, \, D \,], \ \  r(BMB^{*}) + r(D) ]\,  \right\},
\label{}
\\
 \min_{X\in {\mathbb C}^{p \times m}}\!\!\!r[\, \phi(X)\,] & =
\max\{\, 2r[\, A, \, D \,]  - r(M), \ \  r(D) - r(BMB^{*}),  \  r[\, A, \, D \,]  + i_{-}(D)  - i_{+}(BMB^{*}) - i_{-}(N), \nb
\\
& \ \ \ \  \ \  \ \ \ \ r[\, A, \, D \,]  + i_{+}(D)  - i_{-}(BMB^{*}) - i_{+}(N) \},
 \label{}
\\
\max_{X\in {\mathbb C}^{p \times m}} \!\!\!i_{\pm}[\, \phi(X)\,] & =
 \min \left\{ i_{\pm}(M),  \ \ i_{\pm}(BMB^{*}) + i_{\pm}(D) \right\},
\label{}
\\
\min_{X\in {\mathbb C}^{p \times m}} \!\!\!i_{\pm}[\, \phi(X)\,]& = \max\left\{\, r[\, A, \, D \,] -
 i_{\mp}(N), \ \  i_{\pm}(D) - i_{\mp}(BMB^{*}) \right\}.
 \label{}
\end{align}
\end{corollary}

\begin{corollary}\label{T43b}
Let
\begin{align}
\phi(X) = \left(\, AX + C\,\right)\!M\!\left(\, AX + C \right)^{*} + D,
\label{yy432}
\end{align}
where $A \in {\mathbb C}^{n\times p},$ $C \in {\mathbb C}^{n\times m},$
$D \in {\mathbb C}_{{\rm H}}^{n}$ and $M \in {\mathbb C}_{{\rm H}}^{m}$  are given$,$ and $X\in {\mathbb C}^{p\times m}$
is a variable matrix$,$ and define
\begin{align}
& N_1 =\left[\!\! \begin{array}{ccc} D  + CMC^{*} & A \\ A^{*} & 0
 \end{array} \!\!\right]\!,  \ \ N_2 = \left[\!\! \begin{array}{cccc} D  & CM  & A
 \\ A^{*} & 0 & 0 \end{array} \!\!\right]\!.
\label{yy433}
\end{align}
Then$,$
\begin{align}
\max_{X\in {\mathbb C}^{p \times m}} \!\!\!r[\,\phi(X)\,]  & =  \min  \left\{ \, r[\, A , \, CM, \, D \,], \ \
r(N_1),  \ \ r(M) + r(D)\, \right\},
\label{yy434}
\\
\min_{X\in {\mathbb C}^{p \times m}} \!\!\!r[\,\phi(X)\,] & = 2r[\, A, \, CM, \, D \,]  +
\max\{\,  s_1, \ \  s_2,  \ \ s_3, \ \  s_4 \,\},
\label{yy435}
\\
\max_{X\in {\mathbb C}^{p \times m}} \!\!\!i_{\pm}[\,\phi(X)\,] &  =
 \min \left\{ \, i_{\pm}(N_1),  \ \ i_{\pm}(M) + i_{\pm}(D) \, \right\},
\label{yy436}
\\
\min_{X\in {\mathbb C}^{p \times m}} \!\!\!i_{\pm}[\,\phi(X)\,] & = r[\, D, \, CM, \, A \,]  +
\max\left\{\, i_{\pm}(N_1) - r(N_2), \ \   i_{\pm}(D) - r[\, A, \, D\,]  - i_{\mp}(M)\, \right\},
 \label{yy437}
\end{align}
where
\begin{align*}
& s_1 = r(N_1) - 2r(N_2), \ \  s_2 =  r(D)- 2r[\, A, \ D\,]  - r(M),
\\
& s_3 = i_{+}(N_1) - r(N_2) - i_{-}(D)  - r[\, A, \ D\,]  - i_{+}(M), \ \
 s_4 = i_{-}(N_1) - r(N_2) - i_{+}(D) - r[\, A, \ D\,]  - i_{-}(M).
\end{align*}
\end{corollary}

\begin{corollary}\label{T43c}
Let
\begin{align}
\phi(X) = \left(\, XB + C\,\right)\!M\!\left(\, XB + C \right)^{*} + D,
\label{yy438}
\end{align}
where $B \in {\mathbb C}^{p\times m},$ $C \in {\mathbb C}^{n\times m},$
$D \in {\mathbb C}_{{\rm H}}^{n}$ and $M \in {\mathbb C}_{{\rm H}}^{m}$  are given$,$ and $X\in {\mathbb C}^{n\times p}$
is a variable matrix$,$ and define
\begin{align}
N_1 = \left[\!\! \begin{array}{ccc} D  + CMC^{*} & CMB^{*} \\
 BMC^{*} & BMB^{*} \end{array} \!\!\right]\!, \ \  N_2 = [\, \,BMB^{*}, \  BMC^{*}\,].
\label{439}
\end{align}
Then$,$
\begin{align}
\max_{X\in {\mathbb C}^{p \times m}} \!\!\!r[\,\phi(X)\,]  & = n,
\label{440}
\\
\min_{X\in {\mathbb C}^{p \times m}} \!\!\!r[\,\phi(X)\,] & = \max\{\, 0, \ \  r(N_1) - 2(N_2),  \ \
 i_{+}(N_1) - r(N_2), \ \ i_{-}(N_1) - r(N_2) \,\},
\label{441}
\\
\max_{X\in {\mathbb C}^{p \times m}} \!\!\!i_{\pm}[\,\phi(X)\,] &  =
 \min \left\{ \, n, \ \ i_{\pm}(N_1) \, \right\},
\label{442}
\\
\min_{X\in {\mathbb C}^{p \times m}} \!\!\!i_{\pm}[\,\phi(X)\,] & = n  +
\max\left\{\, 0, \ \  i_{\pm}(N_1) - r(N_2) \, \right\}.
 \label{443}
\end{align}
\end{corollary}

We next solve the two quadratic optimization problems in (\ref{13a}), where the two matrices $\phi(\widehat{X})$
and $\phi(\widetilde{X})$, when they exist, are called the global maximal and minimal matrices of $\phi(X)$ in (\ref{11})
  in the L\"owner partial ordering, respectively.

\begin{corollary}\label{T44}
Let $\phi(X)$ be as given  in {\rm (\ref{11})}$,$  and let $N = \left[\!\! \begin{array}{ccc} D + CMC^{*} & CMB \\  B^{*}MC & B^{*}MB \end{array} \!\!\right]\!.$ Then$,$
there exists an $\widehat{X}\in {\mathbb C}^{p\times m}$ such that
 \begin{align}
 \phi(X)\succcurlyeq \phi(\widehat{X})
 \label{zx228}
\end{align}
  holds for all $X\in {\mathbb C}^{p\times m}$ if and only if
\begin{align}
BMB^{*} \succcurlyeq 0, \ \ {\mathscr R}(CMB^{*}) \subseteq {\mathscr R}(A), \ \ {\mathscr R}(BMC^{*})
\subseteq {\mathscr R}(BMB^{*}).
\label{zx229}
\end{align}
In this case$,$ the following hold$.$
\begin{enumerate}
\item[{\rm(a)}]  The matrix $\widehat{X}\in {\mathbb C}^{p\times m}$ satisfying  {\rm (\ref{zx228})}
    is the solution of the linear matrix equation
\begin{align}
A\widehat{X}BMB^{*} + CMB^{*} =0.
\label{zx230}
\end{align}
Correspondingly$,$
\begin{align}
& \widehat{X}  =  - A^{\dag}CMB^{*}(BMB^{*})^{\dag} + F_AV_1 + V_2E_{BMB^{*}},
\label{zx231}
\\
& \phi(\widehat{X}) = D + CMC^{*} - CMB^{*}(BMB^{*})^{\dag}BMC^{*},
\label{zx232}
\\
& \phi(X) - \phi(\widehat{X}) =  (\, AXB + C\,)MB^{*}(BMB^{*})^{\dag}BM(\, AXB + C\,)^{*},
\label{zx233}
\end{align}
where $V_1$ and $V_2$ are arbitrary matrices$.$

\item[{\rm(b)}]  The inertias and ranks of $\phi(\widehat{X})$ and $\phi(X) - \phi(\widehat{X})$ are given by
\begin{align}
&  i_{+}[\,\phi(\widehat{X})\,] = i_{+}(N) - r(BMB^{*}), \  i_{-}[\,\phi(\widehat{X})\,] = i_{-}(N),  \
 r[\,\phi(\widehat{X})\,] = r(N) - r(BMB^{*}),
\label{zx234}
\\
& i_{+}[\,  \phi(X) - \phi(\widehat{X}) \,] = r[\,  \phi(X) - \phi(\widehat{X}) \,]= r(\, AXBMB^{*} + CMB^{*}\,).
\label{zx235}
\end{align}

\item[{\rm (c)}] The matrix $\widehat{X}\in {\mathbb C}^{p\times m}$ satisfying  {\rm (\ref{zx228})}
is unique if and only if
\begin{align}
r(A) = p, \ \  {\mathscr R}(CMB^{*}) \subseteq {\mathscr R}(A), \ \  BMB^{*}\succ 0;
\label{zx236}
\end{align}
under this condition$,$
\begin{align}
& \widehat{X}  =  - A^{\dag}CMB^{*}(BMB^{*})^{-1},
\label{zx237}
\\
& \phi(\widehat{X}) = D + CMC^{*} - CMB^{*}(BMB^{*})^{-1}BMC^{*},
\label{zx238}
\\
& \phi(X) - \phi(\widehat{X}) =  (\, AXB + C\,)MB^{*}(BMB^{*})^{-1}BM(\, AXB + C\,)^{*}.
\label{zx239}
\end{align}

\item[{\rm (d)}] $\widehat{X} = 0$ is a solution of {\rm (\ref{zx228})} if and only if
$BMB^{*} \succcurlyeq 0$ and $CMB^{*} = 0.$ In this case$,$ $\phi(0) = D + CMC^{*}.$

\item[{\rm (e)}] $\widehat{X} = 0$ is the unique solution {\rm (\ref{zx228})} if and only if
$r(A) = p,$ $CMB^{*} =0$ and $BMB^{*}\succ 0.$ In this case$,$  $\phi(0) = D + CMC^{*}.$

\item[{\rm (e)}] There exists an $\widehat{X}\in {\mathbb C}^{p\times m}$ such that
\begin{align}
\phi(X) \succcurlyeq \phi(\widehat{X}) \succcurlyeq 0
 \label{zx240}
\end{align}
holds for all $X\in {\mathbb C}^{p\times m}$ if and only if
\begin{align}
 {\mathscr R}(CMB^{*}) \subseteq {\mathscr R}(A) \ \  and \ \  N \succcurlyeq 0.
\label{zx241}
\end{align}
In this case$,$ the matrix $\widehat{X}\in {\mathbb C}^{p\times m}$ satisfying  {\rm (\ref{zx240})}
is unique if and only if
\begin{align}
r(A) = p, \ \ {\mathscr R}(CMB^{*}) \subseteq {\mathscr R}(A), \ \   BMB^{*} \succ 0, \ \  N \succcurlyeq 0.
\label{zx242}
\end{align}
\end{enumerate}
\end{corollary}

\noindent  {\bf Proof.} \ Let
$$
\psi(X) = \phi(X) - \phi(\widehat{X}) =  \left(\, AXB + C\,\right)\!M\!\left(\, AXB + C \right)^{*} -
\left(\, A\widehat{X}B + C\,\right)\!M\!\left( A\widehat{X}B + C \right)^{*}.
$$
Then, $\phi(X) \succcurlyeq \phi(\widehat{X})$ is equivalent to $\psi(X) \succcurlyeq 0.$
 Under $A \neq 0,$ we see from Corollary 2.2(j) that $\psi(X) \succcurlyeq 0$ holds for all $X \in {\mathbb C}^{p \times m}$
 if and only if
\begin{align}
\left[\!\! \begin{array}{ccc} CMC^{*} - \left( A\widehat{X}B + C\right)\!M\!\left( A\widehat{X}B + C \right)^{*}
& CMB^{*}  \\ BMC^{*} & BMB^{*} \end{array} \!\!\right] \succcurlyeq 0,
\label{zx243}
\end{align}
which, by Lemma \ref{T14}(e), is further equivalent to
\begin{align}
& BMB^{*} \succcurlyeq 0, \ \  \ {\mathscr R}(BMC^{*}) \subseteq {\mathscr R}(BMB^{*}),  \
\label{zx244}
\\
& CMC^{*} - \left( A\widehat{X}B + C\right)\!M\!\left( A\widehat{X}B + C \right)^{*}  -
 CMB^{*} (BMB^{*})^{\dag}BMC^{*} \succcurlyeq 0.
\label{zx245}
\end{align}
In this case, it is easy to verify
\begin{align}
& CMC^{*} - \left(A\widehat{X}B + C\right)\!M\!\left(\, A\widehat{X}B + C \right)^{*} - CMB^{*} (BMB^{*})^{\dag}BMC^{*} \nb
\\
&  =  - (\, A\widehat{X}BMB^{*} +  CMB^{*}\,)(BMB^{*})^{\dag}(\, A\widehat{X}BMB^{*} +  CMB^{*} \,)^{*},
\label{zx246}
\end{align}
 and therefore, the inequality in (\ref{zx245}) is equivalent to
 $ CMB^{*} + A\widehat{X}BMB^{*} =0$.  By Lemma \ref{T15},
this matrix equation is solvable if and only if ${\mathscr R}(CMB^{*}) \subseteq {\mathscr R}(A)$ and
${\mathscr R}(BMC^{*}) \subseteq {\mathscr R}(BMB^{*})$. In this case, the general solution of the equation
is (\ref{zx231}) by Lemma \ref{T15}, and (\ref{zx245}) becomes
$$
CMC^{*} - \left(\, A\widehat{X}B + C\,\right)\!M\!\left(\, A\widehat{X}B + C \right)^{*} - CMB^{*} (BMB^{*})^{\dag}BMC^{*}  =0.
$$
Thus (\ref{zx232}) and (\ref{zx233}) follow. The results in (b)--(f) follow from (a).    \qquad $\Box$

\medskip

The following corollary can be shown similarly.

\begin{corollary}\label{T45}
Let $\phi(X)$ be as given  in {\rm (\ref{11})}$,$ and let  $N = \left[\!\! \begin{array}{ccc} D + CMC^{*} & CMB \\  B^{*}MC & B^{*}MB \end{array} \!\!\right]\!.$ Then$,$  there exists an $\widetilde{X}\in {\mathbb C}^{p\times m}$ such that
\begin{align}
\phi(X) \preccurlyeq \phi(\widetilde{X})
\label{zx247}
\end{align}
 holds for all $X\in {\mathbb C}^{p\times m}$ if and only if
\begin{align}
BMB^{*} \preccurlyeq 0, \ \ {\mathscr R}(CMB^{*}) \subseteq {\mathscr R}(A), \ \ {\mathscr R}(BMC^{*}) \subseteq
{\mathscr R}(BMB^{*}).
\label{zx248}
\end{align}
In this case$,$  the following hold$.$
\begin{enumerate}
\item[{\rm(a)}] The matrix $\widetilde{X}$ satisfying  {\rm (\ref{zx247})} is the solution of the linear matrix equation
\begin{align}
A\widetilde{X}BMB^{*} + CMB^{*} =0.
\label{zx249}
\end{align}
Correspondingly$,$
\begin{align}
& \widetilde{X}  =  - A^{\dag}CMB^{*}(BMB^{*})^{\dag} + F_AV_1 + V_2E_{BMB^{*}},
\label{zx250}
\\
&\phi(\widetilde{X}) = D + CMC^{*} - CMB^{*}(BMB^{*})^{\dag}BMC^{*},
\label{zx251}
\\
&\phi(X)  - \phi(\widetilde{X}) = (\, AXB + C\,)MB^{*}(BMB^{*})^{\dag}BM(\, AXB + C\,)^{*},
\label{zx252}
\end{align}
where $V_1$ and $V_2$ are  arbitrary matrices$.$

\item[{\rm(b)}]  The inertias and ranks of $\phi(\widetilde{X})$ and $\phi(X) - \phi(\widetilde{X})$ are given by
\begin{align}
&  i_{+}[\,\phi(\widetilde{X})\,] = i_{+}(N), \  i_{-}[\,\phi(\widetilde{X})\,] = i_{-}(N) - r(BMB^{*}), \
 r[\,\phi(\widehat{X})\,] = r(N) - r(BMB^{*}),
\label{zx253}
\\
& i_{-}[\,  \phi(X) - \phi(\widetilde{X}) \,] = r[\,  \phi(X) - \phi(\widetilde{X}) \,]= r(\, AXBMB^{*} + CMB^{*}\,).
\label{zx254}
\end{align}

\item[{\rm (c)}] The matrix $\widetilde{X}\in {\mathbb C}^{p\times m}$ satisfying  {\rm (\ref{zx247})}
is unique if and only if
\begin{align}
r(A) = p, \ \  {\mathscr R}(CMB^{*}) \subseteq {\mathscr R}(A), \ \  BMB^{*} \prec 0.
\label{zx255}
\end{align}
In this case$,$
\begin{align}
& \widetilde{X}  =  - A^{\dag}CMB^{*}(BMB^{*})^{-1},
\label{256}
\\
& \phi(\widetilde{X}) = D + CMC^{*} - CMB^{*}(BMB^{*})^{-1}BMC^{*},
\label{257}
\\
& \phi(X) - \phi(\widetilde{X}) =  (\, AXB + C\,)MB^{*}(BMB^{*})^{-1}BM(\, AXB + C\,)^{*}.
\label{258}
\end{align}

\item[{\rm (d)}] $\widetilde{X} = 0$ is a solution of {\rm (\ref{zx247})} if and only if
$BMB^{*} \preccurlyeq 0$ and $CMB^{*} = 0.$ In this case$,$ $\phi(0) = D + CMC^{*}.$

\item[{\rm (e)}] $\widetilde{X} = 0$ is the unique solution of {\rm (\ref{zx247})} if and only if
$r(A) = p,$ $CMB^{*} =0$ and $BMB^{*} \prec 0.$ In this case$,$  $\phi(0) = D + CMC^{*}.$

\item[{\rm (f)}] There exists an $\widetilde{X}\in {\mathbb C}^{p\times m}$ such that
\begin{align}
\phi(X) \preccurlyeq \phi(\widetilde{X}) \preccurlyeq 0
 \label{259}
\end{align}
holds for all $X\in {\mathbb C}^{p\times m}$ if and only if
\begin{align}
 {\mathscr R}(CMB^{*}) \subseteq {\mathscr R}(A) \ \ and \ \  N \preccurlyeq 0.
\label{260}
\end{align}
In this case$,$ the matrix $\widetilde{X}\in {\mathbb C}^{p\times m}$ satisfying  {\rm (\ref{259})}
is unique if and only if
\begin{align}
r(A) = p, \ \ {\mathscr R}(CMB^{*}) \subseteq {\mathscr R}(A), \ \   BMB^{*} \prec 0, \ \  N \preccurlyeq 0.
\label{261}
\end{align}
\end{enumerate}
\end{corollary}

\section{The convexity and concavity of $\phi(X)$ in (\ref{11})}
\renewcommand{\theequation}{\thesection.\arabic{equation}}
\setcounter{section}{5}
\setcounter{equation}{0}

As usual, the matrix-valued function $\phi(X)$ in (\ref{11}) is said to be convex if and only if
\begin{equation}
\phi\left( \frac{1}{2}X_1 + \frac{1}{2}X_2 \right) \preccurlyeq \frac{1}{2}\phi(X_1) + \frac{1}{2}\phi(X_2)
\label{kk46}
\end{equation}
holds for all $X_1, \ X_2 \in {\mathbb C}^{p\times m}$;  said to be concave if and only if
\begin{equation}
\phi\left( \frac{1}{2}X_1 + \frac{1}{2}X_2 \right) \succcurlyeq \frac{1}{2}\phi(X_1) + \frac{1}{2}\phi(X_2)
\label{kk47}
\end{equation}
holds for all $X_1, \ X_2 \in {\mathbb C}^{p\times m}$. It is easy to verify that
\begin{equation}
\phi\left( \frac{1}{2}X_1 + \frac{1}{2}X_2 \right) - \frac{1}{2}\phi(X_1) - \frac{1}{2}\phi(X_2)
= - \frac{1}{4}A(\,X_1 - X_2\,)BMB^{*}(\,X_1 - X_2\,)^{*}A^{*},
\label{kk48}
\end{equation}
which is a special case of (\ref{11}) as well. Applying Theorem  \ref{T41} to (\ref{kk48}), we obtain the following
result.

\begin{theorem}\label{T51}
Let $ \phi(X)$ be as given in {\rm (\ref{11})} with $A \neq 0$ and $BMB^{*} \neq 0.$  Then$,$
\begin{align}
\max_{ X_1 \neq X_2,\,  X_1, \, X_2\in {\mathbb C}^{p \times m}}\!\!\!
r\!\left[\, \phi\!\left( \frac{1}{2}X_1 + \frac{1}{2}X_2 \right) - \frac{1}{2}\phi(X_1) - \frac{1}{2}\phi(X_2)\, \right] & =
\min  \left\{ r(A), \ \  r(BMB^{*}) \right\},
\label{54}
\\
 \min_{ X_1 \neq X_2, \,  X_1, \, X_2\in {\mathbb C}^{p \times m}}\!\!\!
r\!\left[\, \phi\!\left( \frac{1}{2}X_1 + \frac{1}{2}X_2 \right) - \frac{1}{2}\phi(X_1) - \frac{1}{2}\phi(X_2)\, \right] & =
 \left\{\!\! \begin{array}{ll} 1  & \  BMB^{*} \succ 0  \ and \ r(A) =p
 \\
  1  & \  BMB^{*} \prec 0  \ and \ r(A) =p
 \\
 0  &  \  otherwise
\end{array}
\right.,
 \label{55}
\\
\max_{ X_1 \neq X_2, \, X_1, \, X_2\in {\mathbb C}^{p \times m}}\!\!\!
i_{+}\!\left[\, \phi\!\left( \frac{1}{2}X_1 + \frac{1}{2}X_2 \right) - \frac{1}{2}\phi(X_1) - \frac{1}{2}\phi(X_2)\, \right]  & =
 \min\{\, r(A),  \ \ i_{-}(BMB^{*}) \, \},
\label{56}
\\
\max_{ X_1 \neq X_2, \, X_1, \, X_2\in {\mathbb C}^{p \times m}}\!\!\!
i_{-}\!\left[\, \phi\!\left( \frac{1}{2}X_1 + \frac{1}{2}X_2 \right) - \frac{1}{2}\phi(X_1) - \frac{1}{2}\phi(X_2)\, \right]  & =
 \min\{\, r(A),  \ \ i_{+}(BMB^{*}) \, \},
\label{57}
\\
\min_{ X_1 \neq X_2, \ X_1, \, X_2\in {\mathbb C}^{p \times m}}\!\!\!
i_{+}\!\left[\, \phi\!\left( \frac{1}{2}X_1 + \frac{1}{2}X_2 \right) - \frac{1}{2}\phi(X_1) - \frac{1}{2}\phi(X_2)\, \right]&
= \left\{\!\! \begin{array}{ll} 1  & \ \ BMB^{*} \prec 0  \ and \ r(A) =p
 \\
 0  &  \ \ BMB^{*} \nprec 0 \   or  \ r(A) < p
\end{array}
\right. ,
 \label{58}
 \\
\min_{ X_1 \neq X_2, \ X_1, \, X_2\in {\mathbb C}^{p \times m}}\!\!\!
i_{-}\!\left[\, \phi\!\left( \frac{1}{2}X_1 + \frac{1}{2}X_2 \right) - \frac{1}{2}\phi(X_1) - \frac{1}{2}\phi(X_2)\, \right]&
= \left\{\!\! \begin{array}{ll} 1  & \ \  BMB^{*} \succ 0  \ and \ r(A) = p
 \\
 0  &  \   \ BMB^{*} \nsucc 0 \   or  \ r(A) < p
\end{array}
\right. .
 \label{59}
\end{align}
In consequence$,$ the following hold$.$
\begin{enumerate}
\item[{\rm (a)}]  There exist $X_1, \ X_2 \in {\mathbb C}^{p\times m}$ with $X_1 \neq X_2$ such that
$\phi\left( \frac{1}{2}X_1 + \frac{1}{2}X_2 \right) - \frac{1}{2}\phi(X_1)- \frac{1}{2}\phi(X_2)$
is nonsingular if and only if both $r(A) =n$ and $r(BMB^{*}) \succcurlyeq n.$

\item[{\rm (b)}] There  exist $X_1, \ X_2 \in {\mathbb C}^{p\times m}$ with $X_1 \neq X_2$ such that
 $\phi\left( \frac{1}{2}X_1 + \frac{1}{2}X_2 \right) = \frac{1}{2}\phi(X_1)+ \frac{1}{2}\phi(X_2)$
if and only if  $BMB^{*} \nsucc 0$ and $BMB^{*} \nprec 0,$ or $r(A) <p.$

\item[{\rm (c)}] There exist $X_1, \ X_2 \in {\mathbb C}^{p\times m}$ with $X_1 \neq X_2$  such that
$\phi\left( \frac{1}{2}X_1 + \frac{1}{2}X_2 \right) \succ \frac{1}{2}\phi(X_1) + \frac{1}{2}\phi(X_2)$
if and only if both $BMB^{*} \prec 0$ and $r(A) = n.$

\item[{\rm (d)}] There exist $X_1, \ X_2 \in {\mathbb C}^{p\times m}$ with $X_1 \neq X_2$ such that
$\phi\left( \frac{1}{2}X_1 + \frac{1}{2}X_2 \right) \prec \frac{1}{2}\phi(X_1)+ \frac{1}{2}\phi(X_2)$
if and only if  both $BMB^{*} \succ 0$ and $r(A) = n.$

\item[{\rm (e)}] There exist $X_1, \ X_2 \in {\mathbb C}^{p\times m}$ with $X_1 \neq X_2$ such that $\phi\left( \frac{1}{2}X_1 + \frac{1}{2}X_2 \right) \succcurlyeq  \frac{1}{2}\phi(X_1)+ \frac{1}{2}\phi(X_2)$
 if and only if either $BMB^{*} \nsucc 0$  or $r(A) < p.$

\item[{\rm (f)}] There exist $X_1, \ X_2 \in {\mathbb C}^{p\times m}$ with $X_1 \neq X_2$ such that $\phi\left( \frac{1}{2}X_1 +
\frac{1}{2}X_2 \right) \preccurlyeq  \frac{1}{2}\phi(X_1)+ \frac{1}{2}\phi(X_2)$
 if and only if either $BMB^{*} \nprec 0$  or $r(A) < p.$

\item[{\rm (g)}] $\phi\left( \frac{1}{2}X_1 + \frac{1}{2}X_2 \right)  \succcurlyeq \frac{1}{2}\phi(X_1) + \frac{1}{2}\phi(X_2)$
for all $X_1, \ X_2 \in {\mathbb C}^{p\times m}$ with $X_1 \neq X_2$  if and only if $BMB^{*} \preccurlyeq 0$.

\item[{\rm (h)}] $\phi\left( \frac{1}{2}X_1 + \frac{1}{2}X_2 \right)  \preccurlyeq \frac{1}{2}\phi(X_1) + \frac{1}{2}\phi(X_2)$
for all $X_1, \ X_2 \in {\mathbb C}^{p\times m}$ with $X_1 \neq X_2$  if and only if $BMB^{*} \succcurlyeq 0$.

\item[{\rm (i)}] If $\phi(X) $ is a positive semi-definite matrix-valued function$,$ then $\phi(X)$ is convex.

\item[{\rm (j)}] If $\phi(X) $ is a negative  semi-definite matrix-valued function$,$ then $\phi(X)$ is concave.
\end{enumerate}
\end{theorem}

\section{Semi-definiteness of general Hermitian quadratic matrix-valued functions and solutions of the corresponding
partial ordering optimization problems}
\renewcommand{\theequation}{\thesection.\arabic{equation}}
\setcounter{section}{6}
\setcounter{equation}{0}

As an extension of (\ref{11}), we consider the general  quadratic matrix-valued function
\begin{align}
\phi(\,X_1, \ldots, X_k\,) = \left(\, \sum_{i = 1}^{k}A_iX_iB_i + C\,\right)\!M\!\left(\,
\sum_{i = 1}^{k}A_iX_iB_i  + C \right)^{\!*} +D,
\label{zx41}
\end{align}
where $0 \neq A_i \in {\mathbb C}^{n\times p_i},$ $B_i \in {\mathbb C}^{m_i\times q},$  $C \in {\mathbb C}^{n\times q},$
$D \in {\mathbb C}_{{\rm H}}^{n}$ and $M\in {\mathbb C}_{{\rm H}}^{q}$ are given$,$
and $X_i\in {\mathbb C}^{p_i\times m_i}$ is a variable matrix$,$ $i = 1, \ldots, k.$ We treat it as
a combined non-homogeneous linear and quadratic Hermitian matrix-valued function
$\phi = \tau\circ \psi$:
$$
\psi: {\mathbb C}^{p_1\times m_1}\oplus \cdots \oplus {\mathbb C}^{p_k\times m_k}
\rightarrow {\mathbb C}^{n\times q}, \ \ \ \tau: {\mathbb C}^{n\times q} \rightarrow {\mathbb C}_{{\rm H}}^{n}.
$$
This general quadratic function between matrix space includes many ordinary Hermitian quadratic matrix-valued
functions as its special cases. Because more than one variable matrices occur in (\ref{zx41}), we do not know at current
time how to establish analytical formulas for the extremal ranks and inertias of (\ref{zx41}). In this section,
we only consider the following problems:
\begin{enumerate}
\item[{\rm (i)}] establish necessary and sufficient conditions for  $\phi(\, X_1, \ldots, X_k\,) \succcurlyeq 0$
   $(\, \phi(\,X_1, \ldots, X_k\,) \preccurlyeq 0\,)$ to hold for all $X_1,\ldots, X_k$;

\item[{\rm (ii)}] establish necessary and
sufficient conditions for the existence of  $\widehat{X}_1, \ldots, \widehat{X}_k$ and
$\widetilde{X}_1, \ldots, \widetilde{X}_k$ such that
\begin{align}
& \phi(\,X_1, \ldots, X_k\,)  \succcurlyeq \phi(\,\widehat{X}_1, \ldots, \widehat{X}_k\,), \ \ \ \phi(\,X_1, \ldots, X_k\,)
\preccurlyeq \phi(\,\widetilde{X}_1, \ldots, \widetilde{X}_k\,)
\label{zx42}
\end{align}
hold for all $X_1, \ldots, X_k$ in the L\"owner partial ordering, respectively, and give analytical
expressions of $\widehat{X}_1, \ldots, \widehat{X}_k$ and $\widetilde{X}_1, \ldots, \widetilde{X}_k.$
\end{enumerate}

\begin{theorem}\label{T61}
Let $\phi(\,X_1, \ldots, X_k\,)$ be as given in {\rm (\ref{zx41})}$,$  and define $B^{*} = [\, B_1^{*}, \ldots, B_k^{*} \,].$
Also let
\begin{align}
N = \left[\!\! \begin{array}{ccc} D +  CMC^{*} & CMB^{*}  \\ BMC^{*} &
 BMB^{*}  \end{array} \!\!\right]\!.
\label{zx43}
\end{align}
Then$,$ the following hold$.$
\begin{enumerate}
\item[{\rm(a)}] $\phi(X_1, \ldots, X_k) \succcurlyeq 0$ for all $X_1\in {\mathbb C}^{p_1 \times m_1},
\ldots, X_k\in {\mathbb C}^{p_k \times m_k}$  if and only if
$N  \succcurlyeq 0.$

\item[{\rm(b)}] $\phi(\,X_1, \ldots, X_k\,) \preccurlyeq 0$ for all $X_1\in {\mathbb C}^{p_1 \times m_1},
\ldots, X_k\in {\mathbb C}^{p_k \times m_k}$ if and only if $N  \preccurlyeq 0.$

\item[{\rm (c)}]  There exist $\widehat{X}_1, \ldots, \widehat{X}_k$ such that
\begin{align}
 \phi(X_1, \ldots, X_k) \succcurlyeq \phi(\widehat{X}_1, \ldots, \widehat{X}_k)
\label{zx44}
\end{align}
holds for all $X_1\in {\mathbb C}^{p_1 \times m_1}, \ldots, X_k\in {\mathbb C}^{p_k \times m_k}$ if and only if
\begin{align}
BMB^{*} \succcurlyeq 0, \ \ {\mathscr R}( BMC^{*}) \subseteq {\mathscr R}( BMB^{*}).
\label{zx45}
\end{align}
In this case$,$ the matrices $\widehat{X}_1, \ldots, \widehat{X}_k$ are the solutions of the linear matrix equation
\begin{align}
 \sum_{i = 1}^{k}A_i\widehat{X}_iB_iMB^{*} = - CMB^{*}.
\label{zx46}
\end{align}
Correspondingly,
\begin{align}
& \phi(\,\widehat{X}_1, \ldots, \widehat{X}_k\,)  = D + CMC^{*} - CMB^{*}(BMB^{*})^{\dag}BMC^{*},
\label{zx47}
\\
& \phi(X_1, \ldots, X_k)  - \phi(\,\widehat{X}_1, \ldots, \widehat{X}_k\,)  =
\left(\sum_{i = 1}^{k}A_iX_iB_i + C\right)\!MB^{*}(BMB^{*})^{\dag}BM\!\left(\sum_{i = 1}^{k}A_iX_iB_i + C \right)^{\!*}\!.
\label{zx48}
\end{align}

\item[{\rm (d)}]  There exist  $\widetilde{X}_1, \ldots, \widetilde{X}_k$  such that
\begin{align}
 \phi(\,X_1, \ldots, X_k) \preccurlyeq \phi(\,\widetilde{X}_1, \ldots, \widetilde{X}_k\,)
\label{zx49}
\end{align}
holds for all $X_1\in {\mathbb C}^{p_1 \times m_1}, \ldots, X_k\in {\mathbb C}^{p_k \times m_k}$ if and only if
\begin{align}
BMB^{*} \preccurlyeq 0, \ \ {\mathscr R}( BMC^{*}) \subseteq {\mathscr R}(BMB^{*}).
\label{zx410}
\end{align}
In this case$,$ the matrices $\widetilde{X}_1, \ldots, \widetilde{X}_k$ are
 the solutions of the linear matrix equation
\begin{align}
 \sum_{i = 1}^{k}A_i\widetilde{X}_iB_iMB^{*} = - CMB^{*}.
\label{zx411}
\end{align}
Correspondingly,
\begin{align}
& \phi(\,\widetilde{X}_1, \ldots, \widetilde{X}_k) = D + CMC^{*} - CMB^{*}(BMB^{*})^{\dag}BMC^{*},
\label{zx412}
\\
& \phi(\,X_1, \ldots, X_k) - \phi(\widetilde{X}_1, \ldots, \widetilde{X}_k\,) =
\left(\sum_{i = 1}^{k}A_iX_iB_i + C\right)\!MB^{*}(BMB^{*})^{\dag}BM\!\left(\sum_{i = 1}^{k}A_iX_iB_i + C \right)^{\!*}\!.
\label{zx413}
\end{align}
\end{enumerate}
\end{theorem}

\noindent  {\bf Proof.} \
Rewrite (\ref{zx41}) as
\begin{align}
\phi(\,X_1, \ldots, X_k) &= \left(\, A_1X_1B_1 +  \sum_{i = 2}^{k}A_iX_iB_i  + C\,\right)\!M\!
\left(\,A_1X_1B_1 + \sum_{i = 2}^{k}A_iX_iB_i + C \right)^{\!*}  \nb
\\
& \ \ +  D,
\label{zx414}
\end{align}
and applying Corollary \ref{T44} to it, we see that
$\phi(\,X_1, \ldots, X_k) \succcurlyeq 0$ for all $X_1\in {\mathbb C}^{p_1 \times m_1}$ if and only if
\begin{align}
\left[\!\! \begin{array}{ccc} D  & 0 \\ 0 & 0 \end{array} \!\!\right] +
\left[\!\! \begin{array}{c}  \sum_{i = 2}^{k}A_iX_iB_i + C \\ B_1 \end{array} \!\!\right]\!
M\!\left[\!\! \begin{array}{c} \sum_{i = 2}^{k}A_iX_iB_i  + C \\ B_1 \end{array} \!\!\right]^{*} \succcurlyeq 0
\label{zx415}
\end{align}
for all $X_2\in {\mathbb C}^{p_2 \times m_2}, \ldots, X_k\in {\mathbb C}^{p_k \times m_k}$. Note that
\begin{align}
&\left[\!\! \begin{array}{ccc} D  & 0 \\ 0 & 0 \end{array} \!\!\right] +
\left[\!\! \begin{array}{c}  \sum_{i = 2}^{k}A_iX_iB_i + C \\ B_1 \end{array} \!\!\right]\!
M\!\left[\!\! \begin{array}{c}  \sum_{i = 2}^{k}A_iX_iB_i + C \\ B_1 \end{array} \!\!\right]^{*} \nb
\\
& = \left[\!\! \begin{array}{ccc} D  & 0 \\ 0 & 0 \end{array} \!\!\right] +
\left(\sum_{i = 2}^{k} \left[\!\! \begin{array}{c} A_i \\ 0 \end{array} \!\!\right]X_i B_i + \left[\!\! \begin{array}{c}   C  \\ B_1 \end{array} \!\!\right] \right)\!M\!
\left(\sum_{i = 2}^{k} \left[\!\! \begin{array}{c} A_i \\ 0 \end{array} \!\!\right]X_i B_i  + \left[\!\! \begin{array}{c}   C  \\ B_1 \end{array} \!\!\right]  \right)^{\!*}\!.
\label{zx416}
\end{align}
Applying  Corollary \ref{T44}, we see that this matrix is positive semi-definite for
all $X_2\in {\mathbb C}^{p_2 \times m_2}$ if and only if
\begin{align}
\left[\!\! \begin{array}{ccc} D  & 0 & 0 \\ 0 & 0 & 0 \\ 0 & 0 & 0 \end{array} \!\!\right] +
\left[\!\! \begin{array}{c}\sum_{i = 3}^{k} \left[\!\! \begin{array}{c} A_i \\ 0 \end{array} \!\!\right]X_i B_i
 + \left[\!\! \begin{array}{c}   C  \\ B_1 \end{array} \!\!\right] \\ B_2 \end{array} \!\!\right]\!
M\!\left[\!\! \begin{array}{c}
\sum_{i = 3}^{k} \left[\!\! \begin{array}{c} A_i \\ 0 \end{array} \!\!\right]X_i B_i + \left[\!\! \begin{array}{c}   C  \\ B_1 \end{array} \!\!\right] \\ B_2 \end{array} \!\!\right]^{*}
  \succcurlyeq 0
\label{zx417}
\end{align}
for all $X_3\in {\mathbb C}^{p_3 \times m_3}, \ldots, X_k\in {\mathbb C}^{p_k \times m_k}$.
Thus, we obtain by induction that $\phi(X_1, \ldots, X_k) \succcurlyeq 0$ for all $X_1\in {\mathbb C}^{p_1 \times m_1},
\ldots, X_k\in {\mathbb C}^{p_k \times m_k}$ if and only if
$$
\left[\!\! \begin{array}{cccc} D  + CMC^{*} & CMB_1^{*}  & \ldots & CMB_k^{*}  \\ B_1MC^{*} & B_1MB_1^{*}  & \ldots & B_1MB_k^{*} \\
\vdots & \vdots  & \ddots & \vdots \\
 B_kMC^{*}  & B_kMB_1^{*}  & \cdots & B_kMB_k^{*}\end{array} \!\!\right]  =
 \left[\!\!\begin{array}{ccc} D +  CMC^{*} & CMB^{*}  \\ BMC^{*} &  BMB^{*}
  \end{array} \!\!\right]  \succcurlyeq 0,
$$
establishing (a).

Let
\begin{align}
& \rho(X_1, \ldots, X_k)  \nb
\\
& = \phi(X_1, \ldots, X_k) - \phi(\widehat{X}_1, \ldots, \widehat{X}_k) \nb
\\
& = \left( \sum_{i = 1}^{k}A_iX_iB_i + C \right)\!M\!\left(\sum_{i = 1}^{k}A_iX_iB_i  + C \right)^{\!*} -
\left(\, \sum_{i = 1}^{k}A_i\widehat{X}_iB_i + C \right)\!M\!\left(
\sum_{i = 1}^{k}A_i\widehat{X}_iB_i  + C \right)^{\!*}.
\label{zx418}
\end{align}
Then, $\phi(X_1, \ldots, X_k)  \succcurlyeq \phi(\widehat{X}_1, \ldots, \widehat{X}_k)$ for all
$X_1\in {\mathbb C}^{p_1 \times m_1}, \ldots, X_k\in {\mathbb C}^{p_k \times m_k}$  is equivalent to
\begin{align}
\rho(X_1, \ldots, X_k)  \succcurlyeq 0 \ \ {\rm  for \ all} \ \  X_1\in {\mathbb C}^{p_1 \times m_1}, \ldots,
X_k\in {\mathbb C}^{p_k \times m_k}.
\label{zx419}
\end{align}
From (a), (\ref{zx419}) holds if and only if
\begin{align}
\left[\!\! \begin{array}{cccc}  -
\left(\,  \sum_{i = 1}^{k}A_i\widehat{X}_iB_i + C \right)\!M\!\left(
\sum_{i = 1}^{k}A_i\widehat{X}_iB_i +C  \right)^{\!*} + CMC^{*} & CMB^{*}  \\ BMC^{*} & BMB^{*}\end{array} \!\!\right]  \succcurlyeq 0,
\label{zx420}
\end{align}
which, by (\ref{zx243})--(\ref{zx245}), is further equivalent to
\begin{align}
& BMB^{*} \succcurlyeq 0, \ \ {\mathscr R}( BMC^{*}) \subseteq {\mathscr R}( BMB^{*}),
\label{zx421}
\\
& CMC^{*}  - \left( \sum_{i = 1}^{k}A_i\widehat{X}_iB_i + C \right)\!M\!\left(\,
 \sum_{i = 1}^{k}A_i\widehat{X}_iB_i + C \right)^{\!*} -  CMB^{*}(BMB^{*})^{\dag}BMC^{*}  \succcurlyeq 0.
\label{zx422}
\end{align}
In this case,
\begin{align}
& CMC^{*}  - \left(\sum_{i = 1}^{k}A_i\widehat{X}_iB_i + C \right)\!M\!\left( \sum_{i = 1}^{k}A_i\widehat{X}_iB_i +C \right)^{\!*}
-  CMB^{*}(BMB^{*})^{\dag}BMC^{*}  \nb
\\
&  =  - \left( \sum_{i = 1}^{k}A_i\widehat{X}_iB_i  + C \right)MB^{*}(BMB^{*})^{\dag}BM\left( \sum_{i = 1}^{k}A_i\widehat{X}_iB_i
 + C\right)^{\!*}
\label{zx423}
\end{align}
holds, and therefore, (\ref{zx422}) is equivalent to
$CMB^{*} +  \sum_{i = 1}^{k}A_i\widehat{X}_iB_iMB^{*} = 0.$ This is a general
two-sided linear matrix equation involving $k$ unknown matrices. The existence of solutions
of this equation and its general solution  can be derived from the Kronecker product of matrices. The details
are  omitted here.   Result (d) can be shown similarly.  \qquad $\Box$

\medskip

Two consequences of Theorem \ref{T61} are given below.

\begin{corollary}\label{T62}
Let
\begin{align}
\psi(\,X_1, \ldots, X_k\,) = \sum_{i = 1}^{k}\left(\,A_iX_iB_i + C_i\,\right)\!M_i\!\left(\,A_iX_iB_i  + C_i \,\right)^{*} +D,
\label{624}
\end{align}
where $0 \neq A_i \in {\mathbb C}^{n\times p_i},$ $B_i \in {\mathbb C}^{m_i\times q_i},$  $C \in {\mathbb C}^{n\times q_i},$
$D \in {\mathbb C}_{{\rm H}}^{n}$ and $M_i\in {\mathbb C}_{{\rm H}}^{q_i}$ are given$,$
and $X_i\in {\mathbb C}^{p_i\times m_i}$ is a variable matrix$,$ $i = 1, \ldots, k.$
Also define
 $$
 B = {\rm diag}(\,B_1, \ldots, B_k\,),  \ \ C = [\,C_1, \ldots, C_k\,], \ \ M = {\rm diag}(\,M_1, \ldots, M_k\,), \ \
 N = \left[\!\!\begin{array}{ccc} D +  CMC^{*} & CMB^{*}  \\ BMC^{*} &  BMB^{*}
  \end{array} \!\!\right]\!.
 $$
 Then$,$ the following hold$.$
  \begin{enumerate}
\item[{\rm(a)}]
$\psi(\,X_1, \ldots, X_k) \succcurlyeq 0$ for all $X_1\in {\mathbb C}^{p_1 \times m_1},
\ldots, X_k\in {\mathbb C}^{p_k \times m_k}$  if and only if
$N \succcurlyeq 0.$

\item[{\rm(b)}] $\psi(\,X_1, \ldots, X_k) \preccurlyeq 0$ for all $X_1\in {\mathbb C}^{p_1 \times m_1},
\ldots, X_k\in {\mathbb C}^{p_k \times m_k}$ if and only if $N \preccurlyeq 0.$

\item[{\rm (c)}]  There exist $\widehat{X}_1, \ldots, \widehat{X}_k$ such that
\begin{align}
 \psi(\,X_1, \ldots, X_k) \succcurlyeq \psi(\,\widehat{X}_1, \ldots, \widehat{X}_k)
\label{625}
\end{align}
holds for all $X_1\in {\mathbb C}^{p_1 \times m_1}, \ldots, X_k\in {\mathbb C}^{p_k \times m_k}$ if and only if
\begin{align}
B_iM_iB^{*}_i \succcurlyeq 0, \ \ {\mathscr R}(B_iM_iC^{*}_i) \subseteq {\mathscr R}( B_iM_iB^{*}_i), \ \   i = 1, \ldots, k.
\label{626}
\end{align}
In this case$,$ the matrices $\widehat{X}_1, \ldots, \widehat{X}_k$ satisfying {\rm (\ref{625})} are the solutions of the $k$ linear matrix equations
\begin{align}
A_i\widehat{X}_iB_iM_iB^{*}_i = - C_iM_iB^{*}_i,  \ \  i = 1, \ldots, k.
\label{627}
\end{align}
Correspondingly$,$
\begin{align}
& \psi(\,\widehat{X}_1, \ldots, \widehat{X}_k) = D +  \sum_{i = 1}^{k} C_iM_iC^{*}_i  -
\sum_{i = 1}^{k} C_iM_iB^{*}(B_iM_iB^{*}_i)^{\dag}B_iM_iC^{*}_i,
\label{628}
\\
& \psi(\,X_1, \ldots, X_k) - \psi(\,\widehat{X}_1, \ldots, \widehat{X}_k) =\sum_{i = 1}^{k}\left(\,A_iX_iB_i + C_i\,\right)M_iB^{*}_i(B_iM_iB^{*}_i)^{\dag}B_iM_i\left(\,A_iX_iB_i  + C_i \,\right)^{*}.
\label{629}
\end{align}

\item[{\rm (d)}]  There exist  $\widetilde{X}_1, \ldots, \widetilde{X}_k$  such that
\begin{align}
 \psi(\,X_1, \ldots, X_k) \preccurlyeq \psi(\,\widetilde{X}_1, \ldots, \widetilde{X}_k)
\label{630}
\end{align}
holds for all $X_1\in {\mathbb C}^{p_1 \times m_1}, \ldots, X_k\in {\mathbb C}^{p_k \times m_k}$ if and only if
\begin{align}
B_iM_iB^{*}_i \preccurlyeq 0, \ \ {\mathscr R}(B_iM_iC^{*}_i) \subseteq {\mathscr R}( B_iM_iB^{*}_i), \ \   i = 1, \ldots, k.
\label{631}
\end{align}
In this case$,$ the matrices $\widetilde{X}_1, \ldots, \widetilde{X}_k$  satisfying {\rm (\ref{630})} are
 the solutions of the $k$ linear matrix equations
\begin{align}
A_i\widetilde{X}_iB_iM_iB^{*}_i = - C_iM_iB^{*}_i,  \ \  i = 1, \ldots, k.
\label{632}
\end{align}
Correspondingly$,$
\begin{align}
& \psi(\widetilde{X}_1, \ldots, \widetilde{X}_k) = D +  \sum_{i = 1}^{k} C_iM_iC^{*}_i  -
\sum_{i = 1}^{k} C_iM_iB^{*}(B_iM_iB^{*}_i)^{\dag}B_iM_iC^{*}_i,
\label{633}
\\
& \psi(\,X_1, \ldots, X_k) - \psi(\widetilde{X}_1, \ldots, \widetilde{X}_k) =\sum_{i = 1}^{k}\left(\,A_iX_iB_i + C_i\,\right)M_iB^{*}_i(B_iM_iB^{*}_i)^{\dag}B_iM_i\left(\,A_iX_iB_i  + C_i \,\right)^{*}.
\label{634}
\end{align}
\end{enumerate}
\end{corollary}

\noindent  {\bf Proof.} \
Rewrite (\ref{624}) as
\begin{align}
\psi(\,X_1, \ldots, X_k) & = \left[\, A_1X_1B_1 + C_1, \ldots, A_kX_kB_k + C_k \,\right]\!M
\left[\, A_1X_1B_1 + C_1, \ldots, A_kX_kB_k + C_k \,\right]^{\!*} +D \nb
\\
& = \left[\, A_1X_1[\, B_1, \ldots, 0\,] +  \cdots + A_kX_k[\, 0, \ldots, B_k \,]
 + [\,C_1, \ldots, C_k\,] \,\right]\!M \nb
 \\
 &  \ \ \ \times \left[\, A_1X_1[\, B_1, \ldots, 0\,] +  \cdots + A_kX_k[\, 0,\ldots, B_k \,]
 + [\,C_1, \ldots, C_k\,] \,\right]^{*} +D,
\label{635}
\end{align}
which a special case of (\ref{zx41}). Applying Theorem \ref{T61} to it, we obtain the result desired. \qquad $\Box$

\begin{corollary}\label{T63}
Let
\begin{align}
\psi(\,X_1, \ldots, X_k) &=  \left[\, A_1X_1B_1 + C_1, \ldots, A_kX_kB_k + C_k \,\right]\!M
\left[\, A_1X_1B_1 + C_1, \ldots, A_kX_kB_k + C_k \,\right]^{\!*}  \nb
\\
& \ \ + D,
\label{zx424}
\end{align}
where $0 \neq A_i \in {\mathbb C}^{n\times p_i},$ $B_i \in {\mathbb C}^{m_i\times q_i},$  $C_i \in {\mathbb C}^{n\times q_i},$
$D \in {\mathbb C}_{{\rm H}}^{n}$ and $M\in {\mathbb C}_{{\rm H}}^{q_1 + \cdots + q_k}$ are given$,$
and $X_i\in {\mathbb C}^{p_i\times m_i}$ is variable matrix$,$ $i = 1, \ldots, k.$ Also define
 $$
 B = {\rm diag}(\,B_1, \ldots, B_k\,) \ \  and  \ \ C = [\,C_1, \ldots, C_k\,], \ \
 N = \left[\!\!\begin{array}{ccc} D +  CMC^{*} & CMB^{*}  \\ BMC^{*} &  BMB^{*}
  \end{array} \!\!\right]\!.
 $$
 Then$,$ the following hold$.$
 \begin{enumerate}
\item[{\rm(a)}]
$\psi(\,X_1, \ldots, X_k) \succcurlyeq 0$ for all $X_1\in {\mathbb C}^{p_1 \times m_1},
\ldots, X_k\in {\mathbb C}^{p_k \times m_k}$  if and only if
$N \succcurlyeq 0.$

\item[{\rm(b)}] $\psi(\,X_1, \ldots, X_k) \preccurlyeq 0$ for all $X_1\in {\mathbb C}^{p_1 \times m_1},
\ldots, X_k\in {\mathbb C}^{p_k \times m_k}$ if and only if $N \preccurlyeq 0.$

\item[{\rm (c)}]  There exist $\widehat{X}_1, \ldots, \widehat{X}_k$ such that
\begin{align}
 \psi(\,X_1, \ldots, X_k) \succcurlyeq \psi(\,\widehat{X}_1, \ldots, \widehat{X}_k)
\label{zx425}
\end{align}
holds for all $X_1\in {\mathbb C}^{p_1 \times m_1}, \ldots, X_k\in {\mathbb C}^{p_k \times m_k}$ if and only if
\begin{align}
BMB^{*} \succcurlyeq 0, \ \ {\mathscr R}( BMC^{*}) \subseteq {\mathscr R}( BMB^{*}).
\label{zx426}
\end{align}
In this case$,$ the matrices $\widehat{X}_1, \ldots, \widehat{X}_k$ satisfying {\rm (\ref{zx425})} are the solutions of the linear matrix equation
\begin{align}
 \sum_{i = 1}^{k}A_i\widehat{X}_iB_iMB^{*} = - CMB^{*}.
\label{zx427}
\end{align}
Correspondingly$,$
\begin{align}
& \psi(\,\widehat{X}_1, \ldots, \widehat{X}_k) = D + CMC^{*} - CMB^{*}(BMB^{*})^{\dag}BMC^{*},
\label{zx428}
\\
& \psi(\,X_1, \ldots, X_k) - \psi(\,\widehat{X}_1, \ldots, \widehat{X}_k) \nb
\\
&  =\left[\, A_1X_1B_1 + C_1, \ldots, A_kX_kB_k + C_k \,\right]MB^{*}(BMB^{*})^{\dag}BM\left[\, A_1X_1B_1 + C_1, \ldots, A_kX_kB_k + C_k \,\right]^{*}.
\label{zx429}
\end{align}

\item[{\rm (d)}]  There exist  $\widetilde{X}_1, \ldots, \widetilde{X}_k$  such that
\begin{align}
 \psi(\,X_1, \ldots, X_k) \preccurlyeq \psi(\,\widetilde{X}_1, \ldots, \widetilde{X}_k)
\label{zx430}
\end{align}
holds for all $X_1\in {\mathbb C}^{p_1 \times m_1}, \ldots, X_k\in {\mathbb C}^{p_k \times m_k}$ if and only if
\begin{align}
BMB^{*} \preccurlyeq 0, \ \ {\mathscr R}( BMC^{*}) \subseteq {\mathscr R}( BMB^{*}).
\label{zx431}
\end{align}
In this case$,$ the matrices $\widetilde{X}_1, \ldots, \widetilde{X}_k$  satisfying {\rm (\ref{zx430})} are
 the solutions of the linear matrix equation
\begin{align}
 \sum_{i = 1}^{k}A_i\widetilde{X}_iB_iMB^{*} = - CMB^{*}.
\label{zx432}
\end{align}
Correspondingly$,$
\begin{align}
& \psi(\widetilde{X}_1, \ldots, \widetilde{X}_k) = D + CMC^{*} - CMB^{*}(BMB^{*})^{\dag}BMC^{*},
\label{zx433}
\\
& \psi(\,X_1, \ldots, X_k) - \psi(\widetilde{X}_1, \ldots, \widetilde{X}_k) \nb
\\
&  =\left[\, A_1X_1B_1 + C_1, \ldots, A_kX_kB_k + C_k \,\right]MB^{*}(BMB^{*})^{\dag}BM\left[\, A_1X_1B_1 + C_1, \ldots, A_kX_kB_k + C_k \,\right]^{*}.
\label{zx434}
\end{align}
\end{enumerate}
\end{corollary}

\noindent  {\bf Proof.} \
Rewrite (\ref{zx424}) as
\begin{align}
\psi(\,X_1, \ldots, X_k) & = \left[\, A_1X_1B_1 + C_1, \ldots, A_kX_kB_k + C_k \,\right]\!M
\left[\, A_1X_1B_1 + C_1, \ldots, A_kX_kB_k + C_k \,\right]^{\!*} +D \nb
\\
& = \left[\, A_1X_1[\, B_1, \ldots, 0\,] +  \cdots + A_kX_k[\, 0, \ldots, B_k \,]
 + [\,C_1, \ldots, C_k\,] \,\right]\!M \nb
 \\
 &  \ \ \ \times \left[\, A_1X_1[\, B_1, \ldots, 0\,] +  \cdots + A_kX_k[\, 0,\ldots, B_k \,]
 + [\,C_1, \ldots, C_k\,] \,\right]^{*} +D,
\label{zx435}
\end{align}
which a special case of (\ref{zx41}). Applying Theorem \ref{T61} to it, we obtain the result desired. \qquad $\Box$

Many consequences can be derived from the results in this section. For instance,
\begin{enumerate}
\item[{\rm(i)}] the semi-definiteness and the global extremal matrices in the L\"owner partial ordering of the
following constrained QHMF
\begin{align*}
\phi(X) = \left(\, AXB + C\,\right)\!M\!\left(\, AXB + C \right)^{*} + D  \ \ {\rm s.t.}  \ \ PXQ = R
\end{align*}
can be derived;

\item[{\rm(ii)}] the semi-definiteness and the global extremal matrices in the L\"owner partial
ordering of the following matrix expressions that involve partially specified matrices
$$
\left[\!\! \begin{array}{ccc} A  & B  \\ C & ?  \end{array}
\!\!\right]\!M\!\left[\!\! \begin{array}{ccc} A  & B  \\ C & ?  \end{array} \!\!\right]^{*} +N, \ \  \left[\!\! \begin{array}{ccc} ?  & B  \\ C & ?  \end{array}
\!\!\right]\!M\!\left[\!\! \begin{array}{ccc} ?  & B  \\ C & ?  \end{array} \!\!\right]^{*} +N,  \ \
\left[\!\! \begin{array}{ccc} A  &  ? \\ ? & ?  \end{array}
\!\!\right]\!M\!\left[\!\! \begin{array}{ccc} A  & ?  \\ ? & ?  \end{array} \!\!\right]^{*} +N
$$
can be derived. In particular, necessary and sufficient conditions can be derived for the following inequalities
$$
\left[\!\! \begin{array}{ccc} A  & B  \\ C & ?  \end{array}
\!\!\right]\!\left[\!\! \begin{array}{ccc} A  & B  \\ C & ?  \end{array} \!\!\right]^{*} \preccurlyeq I, \ \  \left[\!\! \begin{array}{ccc} ?  & B  \\ C & ?  \end{array}
\!\!\right]\!\left[\!\! \begin{array}{ccc} ?  & B  \\ C & ?  \end{array} \!\!\right]^{*} \preccurlyeq I,  \ \
\left[\!\! \begin{array}{ccc} A  &  ? \\ ? & ?  \end{array}
\!\!\right]\!\left[\!\! \begin{array}{ccc} A  & ?  \\ ? & ?  \end{array} \!\!\right]^{*}\preccurlyeq I
$$
to always hold in the L\"owner partial ordering.
\end{enumerate}


\section{Some optimization problems on the matrix equation $AXB =C$}
\renewcommand{\theequation}{\thesection.\arabic{equation}}
\setcounter{section}{7}
\setcounter{equation}{0}


Consider the following linear matrix equation
\begin{equation}
AXB = C,
\label{zx51}
\end{equation}
where $A \in \mathbb C^{m \times n}, \ B \in \mathbb C^{p \times q}$ and
 $C \in \mathbb C^{m \times q}$  are given,  and $X \in \mathbb C^{n \times p}$
is an unknown matrix. Eq.\,(\ref{zx51}) is one of the best known
matrix equations in matrix theory.  Many papers on this equation and
its applications can be found in the literature. In the Penrose's
seminal paper \cite{Pen}, the consistency conditions and the general
solution of (\ref{zx51}) were completely derived by using generalized
inverse of matrices. If (\ref{zx51}) is not consistent, people often
need to find its approximation solutions under various optimal
criteria, in particular, the least-squares criterion is ubiquitously
used in optimization problems which almost always admits an explicit
global solution. For (\ref{zx51}), the least-squares solution is
defined to be a matrix $X \in \mathbb C^{n \times p}$ that minimizes
the quadratic objective function
\begin{equation}
{\rm trace}\,[\,(\,C - AXB\,\,)(\,C - AXB\,\,)^{*} \,] = {\rm trace}\,[\,(\,C - AXB\,\,)^{*}(\,C - AXB\,\,) \,].
\label{zx52}
\end{equation}
The normal equation corresponding to  (\ref{zx52}) is given by
\begin{equation}
 A^{*}AXBB^{*} = A^{*}CB^{*},
\label{zx53}
\end{equation}
which is always consistent, and the following result is well known.

\begin{lemma}\label{T71}
The general least-squares solution of ${\rm (\ref{zx51})}$ can be written as
\begin{equation}
X = A^{\dag}CB^{\dag} + F_AV_1 + V_2E_B,
\label{zx54}
\end{equation}
where $V_1, \  V_2\in \mathbb C^{n \times p}$ are  arbitrary$.$
\end{lemma}

Define the two QHMFs in (\ref{zx52}) as
\begin{align}
 \phi_1(X): = (\,C - AXB\,\,)(\,C - AXB\,\,)^{*}, \ \ \  \phi_2(X): = (\,C - AXB\,\,)^{*}(\,C - AXB\,\,).
\label{zx55}
\end{align}
Note that
\begin{align}
r[\phi_1(X)]  = r[\phi_2(X)]  = r(\,C - AXB\,\,).
\label{zx56}
\end{align}
Hence, we first obtain the following result from Lemma \ref{T15a}.

\begin{theorem}\label{TU72}
Let $\phi_1(X)$  and $\phi_2(X)$ be as given in {\rm (\ref{zx55})}$.$ Then$,$
\begin{align}
& \max_{X \in \mathbb C^{n \times p}}\!\!\!r[\phi_1(X)] = \max_{X \in \mathbb C^{n \times p}}\!\!\!r[\phi_2(X)] =
\max_{X \in \mathbb C^{n \times p}}\!\!\!r(\,C - AXB\,\,) =  \min \left\{ r[\, A, \, C \,],  \ \
  r \!\left[\!\! \begin{array}{c} B \\ C \end{array} \!\!\right] \right\},
\label{zx57}
\\
& \min_{X \in \mathbb C^{n \times p}}\!\!\!r[
\phi_1(X)] = \min_{X \in \mathbb C^{n \times p}}\!\!\!r[\phi_2(X)] =
\min_{X \in \mathbb C^{n \times p}}\!\!\!r(\,C - AXB\,\,) =  r[\, A, \, C \,] +
 r \!\left[\!\! \begin{array}{c} B \\ C \end{array} \!\!\right] -
 r \!\left[\!\! \begin{array}{cc} C & A \\ B & 0 \end{array} \!\!\right]\!.
\label{zx58}
\end{align}
\end{theorem}

Applying Theorem \ref{T41} to (\ref{zx55}),  we obtain the following result.

\begin{theorem}\label{T52}
Let $\phi_1(X)$  and $\phi_2(X)$ be as given in {\rm (\ref{zx55})}$.$ Then$,$ the following hold.
\begin{enumerate}
\item[{\rm (a)}] There exists an $\widehat{X}\in {\mathbb C}^{n\times p}$ such that
 $\phi_1(X)\succcurlyeq \phi_1(\widehat{X})$ holds for all $X\in {\mathbb C}^{n\times p}$ if and only if
\begin{align}
{\mathscr R}(CB^{*}) \subseteq {\mathscr R}(A).
\label{zx59}
\end{align}
In this case$,$
\begin{align}
\widehat{X} = A^{\dag}CB^{\dag} + F_AV_1 + V_2E_B, \ \ \
 \phi_1(\widehat{X}) = CC^{*} - CB^{\dag}BC^{*},
\label{zx510}
\end{align}
where $V_1, \, V_2 \in \mathbb C^{n \times p}$ are  arbitrary$.$

\item[{\rm (b)}] There exists an $\widehat{X}\in {\mathbb C}^{n\times p}$ such that
 $\phi_2(X)\succcurlyeq \phi_2(\widehat{X})$ holds for all $X\in {\mathbb C}^{n\times p}$ if and only if
\begin{align}
{\mathscr R}(C^{*}A) \subseteq {\mathscr R}(B^{*}).
\label{zx511}
\end{align}
In this case$,$
\begin{align}
\widehat{X} =  A^{\dag}CB^{\dag} + F_AV_1 + V_2E_B, \ \ \phi_2(\widehat{X}) = C^{*}C - C^{*}AA^{\dag}C,
\label{zx512}
\end{align}
where $V_1, \, V_2 \in \mathbb C^{n \times p}$ are  arbitrary$.$
\end{enumerate}
\end{theorem}

Theorem \ref{T52} also motivates us to obtain the following consequence.

\begin{theorem}\label{T53}
Let $A \in \mathbb C^{m \times n}, \ B \in \mathbb C^{p \times q}$ and
 $C \in \mathbb C^{m \times q}$  be  given$.$ Then$,$ there always exist an $X\in \mathbb C^{n \times p}$
  that satisfies
\begin{align}
& \min_{\succcurlyeq}\left\{\, A^{*}(\,C - AXB\,\,)(\,C - AXB\,\,)^{*} A\, : \, X\in \mathbb C^{n \times p} \right\},
\label{zx513}
\\
& \min_{\succcurlyeq}\left\{\,  B(\,C - AXB\,\,)^{*}(\,C - AXB\,\,)B^{*} \, :  \, X\in \mathbb C^{n \times p} \right\},
\label{zx514}
\end{align}
and the general solution is given by
\begin{align}
& {\rm arg}\!\min_{\succcurlyeq}\left\{\,  A^{*}(\,C - AXB\,\,)(\,C - AXB\,\,)^{*}A \, :  \, X\in \mathbb C^{n \times p}
\,\right\} \nb
\\
& = {\rm arg}\!\min_{\succcurlyeq}\left\{\,  B(\,C - AXB\,\,)^{*}(\,C - AXB\,\,)B^{*} \, :  \, X\in \mathbb C^{n \times p}
 \,\right\}   \nb
 \\
& = {\rm arg}\!\!\min_{\!\!\!\!\!\!\!\!\!X\in \mathbb C^{n \times p}}
 {\rm tr}\,\left[\,(\,C - AXB\,\,)(\,C - AXB\,\,)^{*} \,\right]  \nb
 \\
&  =  A^{\dag}CB^{\dag} + F_AV_1 + V_2E_B,
\label{zx515}
\end{align}
where $V_1$ and $V_2$ are  arbitrary matrices$,$ namely$,$ the
solutions of the three minimization problems in {\rm (\ref{zx515})}
are the same.
\end{theorem}

For (\ref{zx51}), the weighted least-squares solutions with respect to positive  semi-define matrices $M$ and $N$
are defined to be matrices $X \in \mathbb C^{n \times p}$ that satisfy
\begin{equation}
{\rm trace}\,[\,(\,C - AXB\,\,)M(\,C - AXB\,\,)^{*} \,] =\min, \ \ {\rm trace}\,[\,(\,C - AXB\,\,)^{*}N(\,C - AXB\,\,) \,] =\min,
\label{zx516}
\end{equation}
respectively. In this case, the two QHMFs in (\ref{zx516}) are
\begin{align}
\phi_1(X)  = (\,C - AXB\,\,)M(\,C - AXB\,\,)^{*},  \ \  \phi_2(X) = (\,C - AXB\,\,)^{*}N(\,C - AXB\,\,),
\label{zx517}
\end{align}
so that the theory on the ranks and inertias of $\phi_1(X)$  and  $\phi_2(X)$ can be established routinely.

Recall that the least-squares solution of a linear matrix equation is defined by minimizing the trace of certain QHMF. For example, the least-squares solution of the well-known linear
matrix equation
\begin{align}
AXB + CYD = E,
\label{zx518}
\end{align}
where $A, \, B, \, C, \,D$ are given, are two matrices $X$ and $Y$
such that
\begin{align*}
& {\rm trace}\,[\,(\,E - AXB - CYD \,)(\, E - AXB - CYD  \,)^{*} \,]
\\
& = {\rm trace}\,[\,(\,E - AXB - CYD \,)^{*}(\,E - AXB - CYD \,\,) \,] = \min.
\end{align*}
Correspondingly,  solutions to the L\"owner partial ordering minimization problems of the two QHMFs
\begin{align*}
(\,E - AXB - CYD \,)(\,E - AXB - CYD \,)^{*}, \ \ (\,E - AXB - CYD \,)^{*}(\,E - AXB - CYD \,)
\end{align*}
 can be derived from Theorem \ref{T61}.

\section{Concluding remarks}
\renewcommand{\theequation}{\thesection.\arabic{equation}}
\setcounter{section}{8}
\setcounter{equation}{0}

We established in this paper a group of explicit formulas for calculating the global maximal and minimal ranks and inertias of (\ref{11}) when $X$ runs over the whole matrix space. By taking
these rank and inertia formulas as quantitative tools, we characterized many algebraic properties
of (\ref{11}), including solvability conditions for some nonlinear matrix equations and inequalities
generated from (\ref{11}), and  analytical solutions to the two well-known
classic optimization problems on the $\phi(X)$ in the L\"owner partial ordering.
The results obtained and the techniques adopted for solving the matrix rank and inertia optimization problems enable us to
make new extensions of some classic results on quadratic forms, quadratic matrix equations and quadratic
matrix inequalities, and to derive many new algebraic properties of nonlinear matrix functions that can hardly
be handled before. As a continuation of this work, we mention some research problems on QHMFs
for further consideration.
\begin{enumerate}
\item[{\rm (i)}] Characterize algebraic and topological properties of generalized Stiefel manifolds
composed by the collections of all matrices satisfying (\ref{23})--(\ref{26}).

\item[{\rm (ii)}]  The difference of (\ref{11}) at two given matrices
$X, \ X + \Delta X\in {\mathbb C}^{p \times m}$
$$
\phi(X + \Delta X) - \phi(X)
$$
is homogenous with respect to $\Delta X$, so that we can add a restriction on its norm,
for instance, $\|\Delta X \| = \sqrt{{\rm tr}[(\Delta X)(\Delta X)^{*}]}  < \delta.$
In this case, establish formulas for calculating the maximal and minimal ranks and inertias of the difference
with respect to $\Delta X \neq 0$, and use them to analyze the behaviors of $\phi(X)$  nearby $X$.
Also note that any matrix $X = (x_{ij})_{p \times m}$ can be decomposed as $X = \sum_{i =1}^{p}\sum_{j =1}^{m}x_{ij}e_{ij}.$
A precise analysis on the difference is to take $\Delta X = \lambda e_{ij}$ and to characterize the behaviors of
 of the difference by using the corresponding rank and inertia formulas.

\item[{\rm (iii)}]  Denote the real and complex parts of (\ref{11}) as  $\phi(X) = \phi_0(X) +
i\phi_1(X)$, where $\phi_0(X)$  and $\phi_1(X)$ two real quadratic matrix-valued functions satisfying
$\phi_0^T(X) = \phi_0(X)$ and $\phi_1^T(X) =-\phi_1(X)$. In this case, establish establish formulas for calculating the maximal and minimal ranks and inertias of  $\phi_0(X)$  and $\phi_1(X)$, and use them to characterize behaviors of  $\phi_0(X)$  and $\phi_1(X)$.

\item[{\rm (iv)}]  Partition $\phi(X)$ in  (\ref{11}) as
$$
\phi(X) = \left[\!\! \begin{array}{cc} \phi_{11}(X)  & \phi_{12}(X)  \\ \phi_{12}^{*}(X)  & \phi_{22}(X)  \end{array} \!\!\right].
$$
In this case, establish formulas for calculating the maximal and minimal ranks and inertias of the submatrices
$\phi_{11}(X)$  and $\phi_{22}(X)$ with respect to $X$, and utilize them to characterize behaviors
of these submatrices.

\item[{\rm (v)}] Most criteria related to vector and matrix optimizations are constructed via
traces of matrices. An optimization theory for (\ref{11}) can also be established by taking the
trace of (\ref{11}) as an objective function. In such a case, it would be of interest to characterize relations between
the two optimization theories for (\ref{11}) derived from trace and  L\"owner partial ordering.

\item[{\rm (vi)}] Establish formulas for calculating
the extremal ranks and inertias of
$$
 \left(\, AXB + C\,\right)\!M\!\left(\, AXB + C \right)^{*} + D
  \ \ {\rm s.t.}  \ \ r(X) \leqslant k,
$$
where $k \preccurlyeq \min\{\, p, \ \ m \,\}$. This rank-constrained matrix-valued function is equivalent to the following biquadratic matrix-valued function
$$
\left(\, AYZB + C\,\right)\!M\!\left(\, AYZB + C \right)^{*} + D, \ \ Y\in \mathbb C^{p \times k}, \ \  Z\in \mathbb C^{k \times m}.
$$
Some previous results on positive semi-definiteness of biquadratic forms can be found in \cite{Cad,Cho3}.

\item[{\rm (vii)}] Establish formulas for calculating
the maximal and minimal ranks and inertias of
$$
 \left(\, AXB + C\,\right)\!M\!\left(\, AXB + C \right)^{*} + D  \ \ {\rm s.t.}  \ \ PX = Q \ {\rm and/or} \ XR = S.
$$
This task could be regarded as extensions of classic equality-constrained quadratic programming problems.

\item[{\rm (viii)}] For two given QHMFs
$$
\phi_i(X) = \left(\, A_iXB_i + C_i\,\right)\!M\!\left(\, A_iXB_i + C_i \right)^{*} + D_i , \ \ i =1, \ 2
$$
of the same size, establish necessary and sufficient conditions for
$\phi_1(X) \equiv \phi_2(X)$ to hold.

\item[{\rm (ix)}] Note that the QHMF in (\ref{11}) is embed into the congruence transformation
for a block Hermitian matrix consisting of the given matrices. This fact prompts us to construct
 some general nonlinear matrix-valued functions that can be embed in congruence transformations for
 block Hermitian matrices, for instance,
\begin{align*}
&\left[\!\!\begin{array}{ccc}  I_{m_1}  & 0 & 0 \\ B_1X_1  & I_{m_2}  & 0 \\  B_2X_2B_1X_1  & B_2X_2  & I_{m_3}
\end{array}\! \!\right]\!\left[\!\! \begin{array}{ccc}  A_{11}  & A_{12} & A_{13} \\ A_{12}^{*}  & A_{22}
& A_{23} \\ A_{13}^{*} & A_{23}^{*} & A_{33} \end{array} \!\right]\!\left[\!\!\begin{array}{ccc} I_{m_1} & X^{*}_1B^{*}_1 &  X_1^{*}B^{*}_1X^{*}_2B^{*}_2
\\ 0 & I_{m_2}  &   X_2^{*}B^{*}_2 \\  0 & 0 & I_{m_3}
\end{array}\!\!\right]
\\
&  = \left[\!\! \begin{array}{ccc}
*  & * & *
\\
*  & *  &  *
\\
* & *  & \phi(\,X_1, \, X_2\,) \end{array} \!\right]\!,
\end{align*}
where
$$
 \phi(\,X_1, X_2\,) = [\, B_2X_2B_1X_1, \,  B_2X_2, \,  I_{m_3} \,]\left[\!\! \begin{array}{ccc}  A_{11}  & A_{12} & A_{13} \\ A_{12}^{*}  & A_{22}
& A_{23} \\ A_{13}^{*} & A_{23}^{*} & A_{33} \end{array} \!\right] \!\left[\!\!\begin{array}{c}  X_1^{*}B^{*}_1X^{*}_2B^{*}_2
\\  X_2^{*}B^{*}_2 \\ I_{m_3}
\end{array}\!\!\right]\!,
$$
which is a special case of the following  nonlinear matrix-valued function
$$
\phi(X_1, \, X_2) = \left(\, A_1X_1B_1 + C_1\,\right)\left(\, A_2X_2B_2 + C_2\,\right)\!M\!\left(\, A_2X_2B_2 + C_2 \right)^{*}\left(\, A_1X_1B_1 + C_1\,\right)^{*} + D.
$$
In these cases, it would be of interest to establish possible formulas for calculating the extremal ranks
and inertias of these nonlinear matrix-valued functions (biquadratic matrix-valued functions), in particular,
to find criteria of identifying semi-definiteness of these nonlinear matrix-valued functions, and to solve
the L\"owner partial ordering optimization problems.

\item[{\rm (x)}]  Two  special forms of (\ref{zx41}) and  (\ref{624}) by setting $X_1 = \cdots = X_k = X$ are
$$
\left(\, \sum_{i = 1}^{k}A_iXB_i + C \,\right)\!M\!\left(\, \sum_{i = 1}^{k}A_iXB_i  +C  \right)^{\!*} +D,
 \ \ \sum_{i = 1}^{k}\left(\,A_iXB_i + C_i\,\right)\!M_i\!\left(\,A_iXB_i  + C_i \,\right)^{*} +D.
$$
In this case, find criteria for the QHMF to be semi-definite, and solve for its
 global extremal matrices in the L\"owner partial ordering.

\item[{\rm (xi)}] Many expressions that involve matrices and their generalized inverses can be
represented as quadratic matrix-valued functions, for instance,
$$
D - B^{*}A^{\sim}_r B, \ \ A- BB^{-}A(BB^{-})^{*}, \ \ A - BB^{-}A - A(BB^{-})^{*} +  BB^{-}A(BB^{-})^{*}.
$$
In these cases, it would be of interest to establish formulas for calculating the maximal and minimal ranks
and inertias of these matrix expressions with respect to the reflexive Hermitian $g$-inverse $A^{\sim}_r$¡¡
of a Hermitian matrix $A$, and  $g$-inverse $B^{-}$ of $B$. Some recent work on the ranks and inertias of the Hermitian
Schur complement $D - B^{*}A^{\sim}B$ and their applications was given in \cite{LT-jota,T-ela12}.
\end{enumerate}

Another type of subsequent work is to reasonably extend the results in the precious sections
 to the corresponding operator-valued functions, for which less quantitative methods are allowed to use.\\

\end{document}